\documentclass[12pt]{amsart}
\usepackage{amsmath}\usepackage[nobysame]{amsrefs}
\usepackage{amsfonts}\usepackage{verbatim}
\usepackage{amssymb}\usepackage{hyperref}
\usepackage{marvosym}\usepackage{wasysym}
\usepackage{dictsym}
\usepackage{amstext}\usepackage{bbold}
\usepackage{amsbsy}
\usepackage{amsopn}\usepackage{MnSymbol}
\usepackage{amsthm}\usepackage[usenames,dvipsnames]{color}

\theoremstyle{definition}

\def\proclaim#1{\vskip0.5em\noindent{\bf #1}\it}
\def\endproclaim{\vskip0.5em\par\noindent\rm}
\def\proclaim#1{\vskip0.5em\noindent{\bf #1}\it}
\def\endproclaim{\vskip0.5em\par\noindent\rm}
\def\demo#1{\vskip0.5em\noindent{\bf #1\ }}

\def\text#1{\mbox{#1}}
\def\flushpar{\par\noindent}

\def\tag#1{\eqno{(#1)}}
\def\mod{\mbox{ mod }}
\newcommand{\mapright}[1]{%
    \smash{\mathop{%
        \hbox to 1cm{\rightarrowfill}
        }
    \limits^{#1}
    }
}
\newcommand{\mapleft}[1]{%
    \smash{\mathop{%
        \hbox to 1cm{\rightarrowfill} $\frak w$ $\D$-nomalized 
        }
    \limits_{#1}
    }
}

\def\e{\epsilon}
\def\a{\alpha}
\def\b{\beta}
\def\G{\Gamma}
\def\g{\gamma}
\def\d{\delta}
\def\D{\Delta}

\def\Si{\Sigma}

\def\l{\lambda}
\def\x{\times}

\def\f{\flushpar}

\def\om{\omega}
\def\Om{\Omega}
\def\B{\mathcal B}

\def\({\left(}
\def\){\right)}

\def\<{\langle}
\def\>{\rangle}
\def\bul{\smallskip\f$\bullet\ \ \ $}
\def\Smi{\smallskip\f{ \Smiley\ \ \ }}
\def\lfl{\lfloor}\def\rfl{\rfloor}\def\sms{\smallskip\f}\def\Par{\smallskip\f\P}
\def\pf{\smallskip\f{\it Proof}\ \ \ \ }\def\sms{\smallskip\f}\def\sms{\smallskip\f}\def\Smi{\smallskip\f{\Large\bf\Smiley\ \ \ }}
\def\lra{\longrightarrow}\def\Lra{\Longrightarrow}\def\pprime{\prime\prime}
\def\lfl{\lfloor}\def\rfl{\rfloor}\def\lcl{\lceil}\def\rcl{\rceil}\def\xbm{(X,\B,m)}
\def\xbmt{(X,\B,m,T)}\def\omft{(\Om,\mathcal F,P,T)}
\def\xyr{\xrightarrow}\def\xyl{\xleftarrow}

\begin{document}
\title{Distributional limits of positive, ergodic stationary processes
and infinite ergodic transformations.}
\author{ Jon. Aaronson and Benjamin Weiss}
\address[Aaronson]{School of Math. Sciences, Tel Aviv University,
69978 Tel Aviv, Israel.}
\email{aaro@tau.ac.il}
\address[Weiss]{Institute of Mathematics
    Hebrew Univ. of Jerusalem,
    Jerusalem 91904, Israel}
\email{weiss@math.huji.ac.il}
\begin{abstract}In this note we identify the distributional limits of non-negative, ergodic stationary processes, showing that all are possible. 
  Consequences for infinite ergodic theory are also explored and new examples of distributionally stable- and $\a$-rationally ergodic transformations are presented.
 
\end{abstract}\thanks{\copyright 2016. Aaronson{\tiny$^\prime$}s research was partially supported by ISF grant No. 1599/13.}
\maketitle\markboth{ Distributional limits}{Aaronson and Weiss}
\section*{\S0 Short Introduction}
\

Classical central limit theory is concerned with the  distributional  convergence of normalized partial sums $\tfrac1{a_n}\sum_{k=1}^nX_n$ of independent, 
identically distributed random variables
$(X_1,X_2,\dots)$.
\

Here, we consider this asymptotic distributional  behavior of normalized partial sums $\tfrac1{a_n}\sum_{k=1}^nX_n$ of random variables $(X_1,X_2,\dots)$ 
generated by a {\it stationary process} ({\tt SP}) by which we mean a 
quintuple $(\Om,\mathcal F,P,T,f)$ where $\omft$ is a  probability preserving transformation  ({\tt PPT}) and $f:\Om\to\Bbb R$ is measurable; the ``generated random variables'' 
being the sequence of random variables $(X_n=f\circ T^n)_{n\ge 0}$ defined on the sample space $(\Om,\mathcal F,P)$.
\

The stationary process $(\Om,\mathcal F,P,T,f)$ is {\it non-negative} if $f\ge 0$; and {\it ergodic} ({\tt ESP}) 
if the underlying {\tt PPT} $\omft$ is an  ergodic {\tt PPT} ({\tt EPPT}). 

\

For independent processes, the possible probability distributions (or laws) occurring as limits  were determined by Paul Levy in \cite{Levy}. They are the 
{\tt stable laws}  
(including the normal distribution of the    central limit theorem).
\

For a general {\tt ESP}, it was shown in  \cite{JPT-BW} that any probability distribution on $\Bbb R$ is a possible limit.
\

This paper is about what happens when the stationary process is non-negative.

\

Our main result on stationary processes is  
\proclaim{Theorem 2} \ \ Let $\omft$ be a {\tt EPPT} and let $Y\in\text{\tt RV}\,(\Bbb R_+)$, then $\exists$
$1$-regularly varying, convex  function $b:\Bbb R_+\to\Bbb R_+$ 
and  a positive measurable function $f:\Om\to\Bbb R_+$
\

so that
\begin{align*}
 \tag{\dsrailways}\label{choochoo}\frac1{b(n)}\sum_{k=0}^{n-1}f\circ T^k\xrightarrow[n\to\infty]{\frak d}\ Y.
\end{align*}

\endproclaim
Here and throughout, 
\bul $\Bbb R_+:=(0,\infty)$, 
\bul for a metric space $Z$,\ $\text{\tt RV}(Z)$ denotes the collection of $Z$-valued random variables, and
\bul $\xrightarrow[n\to\infty]{\frak d}$ denotes {\tt strong distributional  convergence} as defined in  \S1 below.

\

Given a random variable, we'll first construct (theorem 1) a specific {\tt ESP}  satisfying {\it inter alia} (\dsrailways). This will be done by {\tt stacking}. We'll then show that  a general {\tt EPPT}  
induces an extension of the given underlying {\tt EPPT} and that this enables transference of (\dsrailways).
\

Previous work on distributional limits of stochastic processes over arbitrary {\tt EPPT}s can be found in \cite{RB-MD},\cite{Volny},\cite{JPT-BW}.
\

We then apply our results to give new examples of distributionally stable {\tt MPT}s (measure preserving transformations). 
\

In theorem 3  we show ({\it inter alia}) that:
  for any $Y\in\text{\tt RV}\,(\Bbb R_+)$,  $\exists$ a {\tt MPT} 
$\xbmt$ and  a $1$-regularly varying function $a:\Bbb R_+\to\Bbb R_+$ satisfying 
$$\frac1{a(n)}\sum_{k=1}^nf\circ T^k\xrightarrow[n\to\infty]{\frak d}\ Y\int_Xfdm\ \ \forall\ f\in L^1(m)_+.$$
A full statement of theorem 3 is given in \S1 below.
\subsection*{Remarks}\ \ 
\par 1)  It is natural to ask what would be the possible limit laws
of the the partial
sums  of nonnegative {\tt ESP} which are scaled and also centered by positive constants. 

\

That is, what are the possible limit laws of
$\frac{S_n-a(n)}{b(n)}$ where $S_n$ is the $n^{\text{\tiny th}}$ partial sum of a nonnegative {\tt ESP}, and $b(n),\ a(n)>0\ \ (n\ge 1)$ are constants?

\

Our result  shows that any probability distribution
with support bounded from below can be obtained in this fashion. It
is likely that our proof can be modified so as to obtain all
distributions as limits
of these normalized and "centered" sums.  We thank  the referee for raising this issue.
\par 2)  It is also natural to ask about the stochastic processes ocurring as distributional limits of
the random step functions $\Phi_n\in D([0,1])$ (as in \cite{Billingsley}, chapter 3) generated by the partials sums of an {\tt ESP} and defined by
$\Phi_n(t):=\frac{S_{[nt]}}{b(n)}.$
\

For example, if $\frac{S_{n}}{b(n)}\xrightarrow[n\to\infty]{\frak d}\ Y$ as in theorem 2, then, due to the
1-regular variation of $b$, 
$\Phi_n\xrightarrow[n\to\infty]{\frak d}\ \mathcal L_Y$ in $D([0,1])$ where $\mathcal L_Y(t):=tY$.
\section*{\S1 Longer introduction}
\subsection*{ Distributional convergence}
 Consider the compact metric space $([0,\infty],\rho)$ with 
 $$\rho(x,y):=|\tan^{-1}(x)-\tan^{-1}(y)|.$$
 For $x,y\in\Bbb R_+,\ \rho(x,y)\le |x-y|$.
 We'll use the 
 \bul  $\rho$-{\it uniform} distance on $\text{\tt RV}(\Bbb R_+)$ defined   by
$$\frak u(Y_1,Y_2):=\min\,\{\sup\rho(Z_1,Z_2):\  Z=(Z_1,Z_2)\in\text{\tt RV}(\Bbb R_+\x \Bbb R_+),\ Z_i\overset{\frak d}=Y_i\ (i=1,2)\};$$
and the 
\bul $\rho$-{\it Vasershtein} distance  on $\text{\tt RV}(\Bbb R_+)$ defined (as in \cite{Vas})  by
 $$\frak v(Y_1,Y_2):=\min\,\{E(\rho(Z_1,Z_2)):\  Z=(Z_1,Z_2)\in\text{\tt RV}(\Bbb R_+\x \Bbb R_+),\ Z_i\overset{\frak d}=Y_i\ (i=1,2)\}.$$
 Evidently $\frak v(Y_1,Y_2)\le \frak u(Y_1,Y_2)$ and, if $\frak v(Y_1,Y_2)< \e$, then 
 $\exists\ Z=(Z_1,Z_2)\in\text{\tt RV}(\Bbb R_+\x \Bbb R_+),\ Z_i\overset{\frak d}=Y_i\ (i=1,2)$ so that
 $$\text{\tt Prob}\,(\rho(Z_1,Z_2)>\sqrt{\e})<\sqrt{\e}.$$
 For $Y_n,\ Y\in\text{\tt RV}\,(\Bbb R_+)$,
 $$E(g(Y_n))\xrightarrow[n\to\infty]{}\ E(g(Y))\ \forall\ g\in C_B(\Bbb R_+)\ \iff\ \frak v(Y_n,Y)\xrightarrow[n\to\infty]{}\ 0.$$
 See the {\tt Skorohod representation theorem} in \cite{Skorohod} and \cite{Billingsley}.
 \subsection*{Strong distributional convergence}
\

For $(X,\B)$ be a measurable space, we denote the  collection of  probability measures on  $(X,\B)$  by $\mathcal P(X,\B)$.
\

Now let $(X,\B,m)$ be a measure space, $Z$ be a metric space, $F_n:X\to Z$ be measurable,
$Y\in\text{\tt RV}\,(Z)$ and 
$P\in\mathcal P(X,\B),\ P\ll m.$ We'll write
$$F_n\overset{P-\mathfrak d}{\underset{n\to\infty}\lra} Y$$ if
$$ \int_Xg(F_n)dP\xrightarrow[n\to\infty]{}\ E(g(Y))\ \forall\ g\in C_B(Z)$$
and say
(as in \cite{distlim}, \cite{IET} and  \cite{TZ}) that $F_n$ {\it converges strongly in distribution}
(written $F_n\overset{\mathfrak d}{\underset{n\to\infty}\lra} Y)$ if
$$F_n\overset{P-\mathfrak d}{\underset{n\to\infty}\lra} Y\ \forall\   \ P\in\mathcal P(X,\B),\ P\ll m.$$
This is called {\it mixing distributional convergence} in \cite{Renyi} and \cite{Eagleson}.

\

In ergodic situations, strong distributional convergence of normal partial sums is an automatic consequence of distributional convergence. Namely:

\proclaim{Eagleson's Theorem\ \ \cite{Eagleson} {\rm (see also \cite{distlim}, \cite{Aldous}\  and \cite{IET}) }}
\

If $(X,\B,m,T,f)$ is an $\Bbb R$-valued, {\tt ESP}, $a(n)\to\infty\ \text{and}\ \exists\ P\in\mathcal P(X,\B)\ P\ll m$ so that
$$ \int_Xg(\tfrac{S_n}{a(n)})dP\xrightarrow[n\to\infty]{}\ E(g(Y))\ \forall\ g\in C([0,\infty])$$
where $S_n:=\sum_{k=1}^nf\circ T^k$,
then $\tfrac{S_n}{a(n)}\xyr[n\to\infty]{\mathfrak d}Y$.\endproclaim
\subsection*{Examples}
\

\Par1\  Let $\gamma\in (0,1]$ and let  $(\Om,\mathcal A,P,S,f)$ be a positive {\tt SP} where $(f\circ S^n:\ n\ge 1)$ are 
independent random variables  satisfying
$$E(f\wedge t) \underset{t\to\infty}\propto\frac{t}{A(t)}$$ where $A(t)$ $\g$-regularly varying in the sense that
$\frac{A(xt)}{A(t)}\xrightarrow[t\to\infty]{}x^\g\ \forall\ x>0$ (see \cite{BGT}).
                                                                                                                 
  By the stable limit theorem ( \cite{Levy}, also   e.g. XIII.6 in \cite{Feller})

\begin{align*}\tag{\tt SLT }&\frac1{A^{-1}(n)}\sum_{k=1}^nf\circ S^k\ \ \xrightarrow[n\to\infty]{\text{\tt dist}}\ \  Z_\g\ \ 
\end{align*}
where $Z_\g$ is  {\it normalized, $\g$-stable} in the sense that $E(e^{-pZ_\g})=e^{-c_\g p^\g}$ where $c_\g>0$ and $E(Z_\g^{-\g})=1$.
 Note that $Z_1\equiv 1$. For generalizations of this to weakly dependent {\tt SP}s, see \cite{AZ} and references therein.
\

\Par2 \ In \cite{AS}  positive {\tt ESP}s $(\Om,\mathcal F,P,R,f)$ were constructed so that
 $$\frac1{b(n)}\sum_{k=0}^{n-1}f\circ R^k\ \xrightarrow[n\to\infty]{\text{\tt dist}}\ \  \ e^{\frac{1}{2}\mathcal N(0,1)^2}$$
 where $b(n)\propto n\sqrt{\log n}$ and $\mathcal N(0,1)$ is standard normal. 
For example   $R=\tau^f$ where $\tau$ is the dyadic adding machine on $\{0,1\}^\Bbb N$ and 
$f(x):=\min\,\{n\ge 1:\ \sum_{k\ge 1}[(\tau^nx)_k-x_k]=0\}$ is the {\it exchangeability waiting time}. 
\

The following is the main  construction enabling theorem 2. It is a specific construction tailored to the target random variable.
\proclaim{Theorem 1}\ \ Let $Y\in\text{\tt RV}\,(\Bbb R_+)$, then $\exists$ 
\bul an odometer $\xbmt$, 
\bul an increasing, $1$-regularly varying  function $b:\Bbb R_+\to\Bbb R_+$ 
\bul a positive measurable function $f:X\to\Bbb R_+$
\

so that
\begin{align*}
 &\tag{\dsrailways}\frac1{b(n)}\sum_{k=0}^{n-1}f\circ T^k\xrightarrow[n\to\infty]{\frak d}\ Y\\ & \exists\ M>1,\ r>0\ \text{and}\ N_0\ge 1\ \text{such that}\\ &
 \tag{\Bicycle }\label{bike}P([\sum_{k=0}^{n-1}f\circ T^k<xb(n)])\le P(Y\le Mx)\ \ \forall\ x\in (0,r),\ n\ge N_0.
\end{align*}
\endproclaim
The (\dsrailways)   \  condition (repeated from page \pageref{choochoo})  is used in the proofs of theorem 2 and 3.  The (\Bicycle)\  condition will be used in theorem 3 in \S6 to obtain 
examples of $\a$-{\tt rational ergodicity}.
\

The next proposition explains why the normalizing constants are necessarily $1$-regularly varying when the support of $Y$ is compact in $\Bbb R_+$.
\proclaim{Normalizing constant proposition}
\

Suppose that $(\Om,\mathcal F,P,R,f)$  is a positive {\tt ESP}, $b(n)>0$, and 
$Y\in\text{\tt RV}\,(\Bbb R_+)$ with $\min\text{\tt supp}\,Y=:a>0\ \text{and}\ \max\text{\tt supp}\,Y=:b<\infty$.
\

If
$\frac{S_n}{b(n)}\xyr[n\to\infty]{\frak d}\ Y$ where $S_n:=\sum_{k=1}^nf\circ T^k$, then $b$ is $1$-regularly varying.\endproclaim\demo{Proof}\ \ 
It suffices to show that $\frac{b(2n)}{b(n)}\xyr[n\to\infty]{}\ 2$. To see this, suppose otherwise, then there exist $\e>0$ and 
a subsequence $K\subset\Bbb N$,  so that
\begin{align*}\tag{\ddag}|\frac{b(2n)}{b(n)}-2|\ge \e\ \ \forall\ \ n\in K. 
\end{align*}

Next, by compactness, there is a further subsequence $K'\subset K$ and a random variable $Z=(Z_1,Z_2)\in\text{\tt RV}\,([0,\infty]^2)$ so that
\begin{align*}
\left(\frac{S_n}{b(n)},\frac{S_n\circ T^n}{b(n)}\right)\xyr[n\to\infty]{\frak d}\ Z.
 \end{align*}

By assumption, we have that $\text{\tt dist}\,Z_i=\text{\tt dist}\,Y\ \ (i=1,2)$. Thus,
$$2a\le Z_1+Z_2\le 2b.$$
Now fix $K^{\pprime}\subset K'$ so that $\frac{b(2n)}{b(n)}\xyr[n\to\infty,\ n\in K^{\pprime}]{}\ c\in [0,\infty]$.
\

By assumption,
\begin{align*}Y&\xyl[n\to\infty,\ n\in K^{\pprime}]{\frak d}\ \frac{S_{2n}}{b(2n)}\\ &=
\frac{b(n)}{b(2n)}\left(\frac{S_n}{b(n)}+\frac{S_n\circ T^n}{b(n)}\right)\\ &\xyr[[n\to\infty,\ n\in K^{\pprime}]{\frak d} c^{-1}(Z_1+Z_2). 
\end{align*}
It follows that $c\in\Bbb R_+$ and that $Z_1+Z_2\overset{\text{\tiny dist}}=cY$.
So on the one hand $\min\text{\tt supp}\,cY=ca\ \text{and}\ \max\text{\tt supp}\,cY=:cb<\infty$ and on the other hand,
$$ca=\min\text{\tt supp}\,(Z_1+Z_2)\ge 2a\ \text{and}\ cb=\max\text{\tt supp}\,(Z_1+Z_2)\le 2b$$
with the conclusion that $c=2$ which contradicts (\ddag).\ \ \Checkedbox

\

\subsection*{Distributional convergence in infinite ergodic theory}
\

Let $\xbmt$ be a conservative, ergodic {\tt MPT} ({\tt CEMPT}) and let $Y\in\text{\tt RV}\,([0,\infty])$. Let  $n_k\uparrow\infty$ be a subsequence  and let $d_k>0$ be constants. 
As in \cite{distlim} and \cite{IET}, we'll write
$$\frac{S_{n_k}^{(T)}}{d_k}\xrightarrow[k\to\infty]{\frak d}\ Y$$
if
$$\frac{S_{n_k}^{(T)}(f)}{d_k}\xrightarrow[k\to\infty]{\frak d}\ Y\int_Xfdm\ \forall\ f\in L^1_+.$$
Call the random variable $Y\in\text{\tt RV}\,([0,\infty])$ appearing a {\it subsequence distributional limit of $T$} and let
$$\mathcal L_T:=\{\text{\tt subsequence distributional limits of $T$}\}.$$

The collection 
$$\{T\in\text{\tt MPT}\,(\Bbb R):\ \mathcal L_T=\text{\tt RV}\,([0,\infty])\}$$
is residual in $\text{\tt MPT}\,(\Bbb R)$, the group of invertible transformations of $\Bbb R$  preserving  Lebesgue measure, 
equipped with the weak topology (see \cite{gendistlim}).
\

 We
call the {\tt CEMPT}
 $(X,\B,m,T)$   {\it  distributionally stable} if there are constants $a(n)=a_{n,Y}(T)>0$  and a random variable $Y$ on $(0,\infty)$ 
 (called the {\it ergodic limit})
 so that
 \begin{align*}
  \tag{\Football}\label{football}\frac{S_n^{(T)}}{a(n)}\xyr[n\to\infty]{\mathfrak d} \ Y.
 \end{align*}
The sequence of constants $(a_{n,Y}(T):\ n\ge 1)$ is determined up to asymptotic equality and we call it the {\it $Y$-distributional return sequence}.
Note that $a_{n,cY}(T)\sim \frac1ca_{n,Y}(T)$.
For distributionally stable {\tt CEMPT}s which are also weakly rationally ergodic, we have that $a_{n,Y}(T)\propto a_n(T)$ the usual return sequence (see  \cite{RE}).
\

Classic examples of distributionally stable {\tt CEMPT}s are obtained via
the Darling-Kac theorem (\cite{DK}): pointwise dual ergodic transformations (e.g. Markov shifts) with 
regularly varying return sequences are distributionally stable with Mittag-Leffler ergodic  limits (see also \cite{IET}, \cite{distlim}).
\

More recently, it has been shown that certain ``random walk adic'' transformations have exponential chi-square distributional limits (see \cite{AS},  \cite{ADDS}  and \cite{Bromberg}).
\

Our main result about infinite, ergodic transformations is
\proclaim{Theorem 3} \ \  For each  $Y\in\text{\tt RV}\,(\Bbb R_+)$,  there is a distributionally stable {\tt CEMPT}  $(X,\B,m,T)$ with ergodic limit $Y$
with $a_{n,Y}(T)$ $1$-regularly varying and 
 $\Om\in\B,\ m(\Om)=1$ so that
\begin{align*}\tag{\dsaeronautical}\label{aero}m(\Om\cap [S_n(1_\Om)\ge xa(n)])\le 2P(Y\text{$\ge$} x)\ \forall\ x>1\ \text{and}\ n\ge 1\ \text{large}.\end{align*}
\endproclaim
The (\dsaeronautical)\  condition (which is an inversion of the (\Bicycle) condition on page \pageref{bike})  will be used in the construction of $\a$-rationally ergodic {\tt MPT}s in \S6.
\

By proposition 3.6.3 in \cite{IET}, distributional stability of a {\tt CEMPT} entails existence of a {\tt law of large numbers} (as in \cite{distlim} and \cite{IET})
for it. An example in \S6 shows it does not entail $\a$-rational ergodicity. 
\subsection*{Plan of the paper}
\

In \S2, we recall the {\tt stacking method} used to construct the odometer in theorem 1. This odometer is constructed together with a sequence of
{\tt step functions} and in \S3, we formulate the {\tt step function extension lemma} needed for the proof of theorem 1 where the limit is a 
{\tt rational random variable}
(taking finitely many values, each with rational probability). In \S4 we prove the step function extension lemma and theorem 1 in this (rational rv) case. In \S5, we prove theorem 1 in general, developing the necessary approximations
of random variables by rational ones. We conclude in \S6 by proving theorem 3 and considering some of its consequences in infinite ergodic theory.
 \section*{\S2 The stacking constructions}

 \
 
{\it Stacking} as in \cite{Chacon} (aka the {\it stacking method} \cite{Friedman} and   {\it cutting and stacking} in \cite{Shields-cutstack},\cite{Shields-book})
is a construction procedure yielding a piecewise translation of an almost open subset  $X\subset \Bbb R$. 
This transformation is invertible and preserves Lebesgue measure.
 \

 As in \cite{Chacon} and \cite{Friedman}, a  {\it column} 
 is a finite sequence of disjoint intervals $W=(I_1,I_2,\dots,I_h)$.
  with equal lengths. The {\it width} of the column is the length of $I_k$. The {\it height} of the column is $h$ and we'll sometimes call
  $W=(I_1,I_2,\dots,I_h)$ an $h$-{\it column}.
  
  \
  
  The {\it base} of the column $W=(I_1,I_2,\dots,I_h)$ is $B(W):=I_1$, its {\it top} is $A(W):=I_h$ and its {\it union} is 
 $U(W)=\bigcupdot_{k=1}^hI_k$.  The {\it measure} of a column is the length of its union. Columns $W\ \text{and}\ W'$ are disjoint if their unions are disjoint. 
  
  \

 The column $W$ is  equipped  with the periodic map 
 $T=T_W:U(W)\to U(W)$ 
 defined by the translations $T:I_k\to I_{k+1}\ \ (1\le k\le h-1)\ \text{and}\ T:I_h\to I_1$.

A {\it castle} ({\it tower} in \cite{Chacon} and \cite{Friedman}) is a finite collection of disjoint columns.

A castle consisting of a single column is known as a {\it   Rokhlin tower}. 
\

A castle is called {\it homogeneous} if all the columns have the same height and width. As before, an homogeneous castle consisting of $h$-columns is called an $h$-{\it castle}.
\

The {\it base}  of the  castle $\frak W=\{W_1,W_2,\dots,W_n\}$ is
$B(\frak W)=\bigcupdot_{k=1}^nB(W_k)$, its top is $A(\frak W)=\bigcupdot_{k=1}^nA(W_k)$  and its union is  $U(\frak W)=\bigcupdot_{k=1}^nU(W_k)$ . 

\

It is equipped with the periodic transformation $T_\frak W:U(\frak W)\to U(\frak W)$ defined by
$$T_{\frak W}|_{U(W_k)}\equiv T_{W_k}.$$
\subsection*{Refinements of castles}
\

The castle $\frak W'$ {\it refines} the castle $\frak W$ (written $\frak W'\succ\frak W$)  if 
\sms (i) each interval of $\frak W$ is a union of intervals of $\frak W'$;
\sms (ii) $A(\frak W')\subset\ A(\frak W)\ \text{and}\ B(\frak W')\subset\ B(\frak W)$;
\sms (iii) $T_{\frak W'}|_{U(\frak W)\setminus A(\frak W)}\equiv T_{\frak W}$.

If $\frak W'\succ\frak W$, then $U(\frak W')\supset U(\frak W)$. 
\

All castle refinements  $\frak W'\succ\frak W$ considered here are  mass preserving in the sense that  $U(\frak W')=U(\frak W)$ (no ``spacers'' are added). 
\

Call the refinement $\frak W'\succ\frak W$ {\it transitive} if 
$$m(U(W')\cap U(W))>0\ \forall\ W'\in \frak W'\ \ \text{and}\ W\in \frak W.$$

\

A sequence $(\frak W_n)_{n\ge 1}$ of castles is a  {\it nested sequence} if each $\frak W_{n+1}$ refines $\frak W_{n}$.

\

Let $(\frak W_n)_{n\ge 1}$ be a nested sequence of castles and consider  the measure space $\xbm$ with $X:=\bigcup_{n=1}^\infty U(\frak W_n)$ equipped with Borel sets $\B$ and Lebesgue measure $m$.

\

As shown in \cite{Chacon} and \cite{Friedman}, 

\Smi There is a  measure preserving transformation  $(X,\B,m,T)$  defined by 
$$T(x)=\lim_{n\to\infty}T_{\frak W_n}(x)\ \ \text{for $m$-a.e.}\ x$$
iff $m(A(\frak W_n))\xrightarrow[n\to\infty]{}\ 0$.

\

It is standard to show that if  infinitely many of the refinements $\frak W_{n+1}\succ\frak W_n$ are transitive, 
then  $(X,\B,m,T)$ is ergodic.

\

The transformation  $(X,\B,m,T)$  is aka the {\it inverse limit} of $(\frak W_n)_{n\ge 1}$ and denoted $T=\varprojlim_{n\to\infty}\frak W_n$.

\

\subsection*{Odometers}
\

An {\it odometer} is  an inverse limit of a (mass preserving) nested sequence of Rokhlin towers. Odometers are ergodic because if $\frak W',\ \frak W$ are  Rokhlin towers and
$\frak W'\ \succ\  \frak W$, then the refinement is clearly transitive. The odometers are the ergodic transformations with  rational, pure point spectrum.
\

\sms{\bf Induced Transformation} (as in \cite{Kakutani})
\

Let $\xbmt$ be a {\tt CEMPT} and let $\Om\in\B,\ 0<m(\Om)<\infty$. The {\it first return time} to $\Om$ is the function  
$\varphi_\Om:\Om\to\Bbb N\cup\{\infty\}$ defined by  $\varphi_\Om(x):=\min\,\{n\ge 1:\ T^nx\in\Om\}$ 
which is finite for a.e. $x\in\Om$ by conservativity.
\

 The {\it induced transformation} is $(\Om,\mathcal B\cap\Om ,m_\Om,T_\Om)$ 
where $T_\Om:\Om\to\Om$ is defined by
$T_\Om(x):=T^{\varphi_\Om(x)}$ and $m_\Om(\cdot):=m(\cdot\|\Om)$. It is  a {\tt PPT}.

\proclaim{Odometer factor proposition}
\

\ \ Let $R$ be an odometer and let $\xbmt$ be an aperiodic {\tt PPT}, then $\exists\ \Om\in\B,\ m(\Om)>0$ so that $R$ is a factor 
{\tt PPT} of $T_\Om$.\endproclaim\demo{Proof}
\

Let $R=\varprojlim_{n\to\infty}\frak W_n$ where $(\frak W_n)_{n\ge 1}$ is a nested sequence of Rokhlin towers. Let the height of $\frak W_n$ be $H_n$, then
there is a sequence $a_1,a_2,\dots\in\Bbb N,\ a_n\ge 2$ so that $H_1=a_1,\ H_{n+1}=a_{n+1}H_n$.\
\

By the basic Rokhlin lemma, for any $\e_1\in (0,1)$ there is  some $B_1$ of positive measure such that the
sets $\{T^i(B_1): i=0,1..a_1-1\}$ are disjoint and
$$X= \bigcupdot_{i=0}^{a_1-1}T^i(B_1)\cupdot E_1$$ where $E_1\in\B$ and $m(E_1)=\e_1m(B_1)$. 
\

Next apply the Rokhlin lemma
again to the induced transformation $T_{B_1}$ with $\e_2\in (0,1)$ to get a base $B_2 \subset
B_1$ with the sets $\{T_{B_1}^iB_2:\ 0\le i<a_2\}$
disjoint and
$$B_1= \bigcupdot_{i=0}^{a_2-1}T_{B_1}^i(B_2)\cupdot E_2$$ where $E_2\in\B(B_1)$  and $m(E_2)=\e_2m(B_2)$.
\

This process is continued to obtain $B_k\in \B,\ B_k\subset B_{k-1}$ with the sets $\{T_{B_{k-1}}^iB_k:\ 0\le i<a_k\}$
disjoint and
$$B_{k-1}= \bigcupdot_{i=0}^{a_k-1}T_{B_{k-1}}^i(B_k)\cupdot E_k$$ where $E_k\in\B(B_{k-1})$  and $m(E_k)=\e_km(B_{k-1})$.
If $\sum_{k\ge 1}\e_k<1$, then
$$\Om:=\bigcap_{k\ge 1}\bigcup_{i=0}^{H_k-1}T^i(B_k)$$
is as advertised.\ \ \Checkedbox

We'll need a condition for an inverse limit of castles to be isomorphic to an odometer.
\

If $W=(I_1,I_2,\dots,I_k)$ and $W'=(I'_1,I'_2,\dots,I'_{k'})$ are  disjoint columns  of intervals with equal width, the {\it stack} of $W\ \text{and}\ W'$ is the column
$$W\ocirc W':=(I_1,I_2,\dots,I_k,I'_1,I'_2,\dots,I'_{k'}).$$

\

Let $q\in\Bbb N$. The column $W$ can be sliced into $q$ subcolumns $$^qW_1,^qW_2,\dots,^qW_q$$ of equal width and the same height.

For a column $W$ and $q\in\Bbb N$, $W^{\circledast q}$ denotes the column obtained from $W$ by slicing the column into
$q$ disjoint subcolumns of equal width and then stacking. That is 
$$W^{\circledast q}=\bigocirc_{k=1}^q {^qW_k}.$$ 

\

Let $\frak W=\{W_k:\ 1\le k\le K\}\ \text{and}\ \ \frak W'=\{W'_\ell:\ 1\le\ell\le L\}$ be  homogeneous castles. 

\

The refinement  $\frak W'\succ\frak W$ is {\it uniform} 
if $\exists\ Q\ge 1,\ \kappa_1,\kappa_2,\dots,\kappa_Q\in\{1,2,\dots,K\}$ with $\{\kappa_q:\ 1\le q\le Q\}=\{1,2,\dots,K\}$ and $s_1,s_2,\dots,s_Q\in\Bbb N$ so that
$$W_\ell'=\  ^L\left(\bigocirc_{q=1}^QW_{\kappa_q}^{\circledast s_q}\right)_\ell.$$

Note that a uniform refinement is  transitive.
\

The nested sequence of homogeneous castles  $(\frak W_n)_{n\ge 1}$ is called {\it uniformly nested} if each refinement  $\frak W_{n+1}\succ\frak W_n$ is uniform. 
\proclaim{Proposition}\ \ Let $(\frak W_n)_{n\ge 1}$ be a uniformly nested sequence of homogeneous castles, then
the {\tt EPPT} $\xbmt:=\varprojlim_{n\to\infty}\frak W_n$ is an odometer.\endproclaim\demo{Proof}
Let $\frak W_n=\{W_j^{(n)}:\ 1\le j\le k_n\}$ and suppose that
$$W_\ell^{(n+1)}=\  ^{k_{n+1}}\left(\bigocirc_{q=1}^{Q_{n+1}}W_{\kappa_q}^{(n)\circledast s^{(n+1)}_q}\right)_\ell,$$
then
$$W_\ell^{(n+1)}=\  ^{k_{n+1}}(\widetilde{W}^{(n)})_\ell$$ where
$$\widetilde{W}^{(n)}:=\bigocirc_{q=1}^{Q_{n+1}}W_{\kappa_q}^{(n)\circledast s^{(n+1)}_q}.$$
The Rokhlin tower $\widetilde{\frak W}^{(n)}:=\{\widetilde{W}^{(n)}\}$ is refined by $\widetilde{\frak W}^{(n+1)}$ and
$$\xbmt=\varprojlim_{n\to\infty}\widetilde{\frak W}^{(n)}.\ \ \CheckedBox$$

\section*{\S3  Step functions, labeled castles and block arrays   }

\

Here we introduce the framework for the proof of Theorem 1.

\

  We'll a construct recursively  a nested sequence of homogeneous, unit measure castles $(\frak W_n)_{n\ge 1}$ and set $\xbmt=\varprojlim_{n\to\infty}\frak W_n$. 
\

The advertised function  $f:X\to\Bbb R_+$ will be defined as $f=\lim_{n\to\infty}f^{(n)}$ where
$f^{(n)}:\frak W_n\to\Bbb R_+$ is a {\it step function} in the sense that it is constant on each of the intervals making up each column in the castle $\frak W_n$.
\

If $\frak W_n=\{W_j^{(n)}:\ 1\le j\le k_n\}$ where each $W_j^{(n)}=(I_{j,k}^{(n)})_{1\le k\le h_n}$ is a column of height $h_n$, then
$$f^{(n)}\cong (w_j^{(n)}:\ 1\le j\le k_n)\subset(\Bbb R_+^{h_n})^{k_n}$$  where 
$$f^{(n)}\equiv w_j^{(n)}(k)\ \text{on}\ \ I_{j,k}^{(n)}.$$

Formally, let 
a $J$-{\it block} be a positive vector $w\in\Bbb R_+^J$  (where $J\in\Bbb N$). The {\it length} of $J$-block $w$ is
$|w|:=J$. 
\

A block  $w\in\Bbb R_+^J$ determines a {\it labeled column}: an {\it underlying column} $W=(I_1,I_2,\dots,I_J)$ together with a step function $F_W:U(W)\to\Bbb R_+$ defined by
$$F_W=\sum_{k=1}^Jw_k1_{I_k}.$$
\

A {\it block array} is an ordered collection of blocks of the same length (called $J$-{\it block array} when all the blocks have length $J$).
\

The block array $\frak w=(w_1,w_2,\dots,w_N)\in\ (\Bbb R_+^h)^N$ determines a  {\it labeled castle}:

\

an {\it underlying castle} $\frak W=(W_1,W_2,\dots,W_N)$  of height $h$, together with a step function $F_\frak w:U(\frak w)\to\Bbb R_+$ defined by
$$F_\frak w:=\sum_{k=1}^N1_{U(W_k)}F_{W_k}.$$

We'll say that the  block array $\frak y$ {\it refines} the block array $\frak x$
written $\frak y\succ \frak x$ if the castle determined by  $\frak y$ refines that determined by $\frak x$.
\

Blocks can be {\tt concatenated}. If $w\in\Bbb R^J\ \text{and}\ w'\in\Bbb R^{J'}$, the {\it concatenation} of $w\ \text{and}\ w'$ is
$$w\odot w':=(w_1,w_2,\dots,w_J,w'_1,w'_2,\dots,w'_{J'})\in\Bbb R^{J+J'}.$$
The concatenation of blocks corresponds to the stacking of their underlying columns.

If $W\ \text{and}\ W'$ are columns of height $J$ and $J'$ respectively and with the same width, and $w\in\Bbb R^J\ \text{and}\ w'\in\Bbb R^{J'}$, then
$$F_{w\odot w'}\equiv F_{\{w, w'\}}\ \ \text{on}\ \ U(W\ocirc W')=U(\{W,W'\})=U(W)\cupdot U(W').$$
\

Similarly, self concatenation $w^{\odot q}$ of the same block $w$ corresponds to cutting and stacking $W^{\ostar q}$ of the corresponding column $W$.

\

\

We call a sequence of block arrays {\it nested} if the underlying sequence of castles is nested.

\

We'll obtain the required {\tt ESP} by producing a nested sequence $(\frak w_n)_{n\ge 1}$ of block arrays whose associated sequence of step functions $(F_{\frak w_n})_{n\ge 1}$
is convergent.

 \subsection*{Block statistics} \

Distributional convergence will be achieved by controlling the empirical distributions of the various short-term partial sums over the tall block arrays.
\

Given a block $w\in\Bbb R_+^h$,  define 
$$S_k(F_w):=\sum_{j=0}^{k-1}F_w\circ T_w^j$$ where $T_w$ is the periodic transformation defined on  the column underlying $w$.
We have
$$S_k(F_w)=\sum_{\nu=1}^hS_k(w)(\nu)1_{I_\nu}$$

where, for $1\le\nu\le h$,
$$S_k(w)(\nu):=\sum_{j=0}^{k-1}w_{\nu+j}.$$
Here translation is considered $\mod h$ that is $\nu+j:=\nu+j\ \mod h$. 

\

For a block array $\frak w=\{w_j:\ 1\le j\le K\}$, set 
$$S_k(F_\frak w)=\sum_{j=1}^K1_{U(w_j)}F_{w_j}$$
and $S_k(\frak w)(\nu,j):=S_k(w_j)(\nu).$
\

We  study the distributions of $S_k(w)\ \text{and}\ S_k(\frak w)$ considered as $\Bbb R_+$-valued random variables on the symmetric probability spaces
$\{1,2,\dots,h\}$ and 
 $\{1,2,\dots,h\}\x \{1,2,\dots,K\}$ respectively.

\

If $w\in\Bbb R^h$ and $m\in\Bbb N$, then
$$S_k(w^{\odot m})(\nu)=S_k(w)(\nu\ \mod h)$$
whence $S_k(w^{\odot m})\ \text{and}\ S_k(w)$ are equidistributed.

In a similar manner, we consider partial sums on a block array $\frak w=\{w_k:\ 1\le k\le n\}:\{1,\dots,h\}\x\{1,\dots,n\}\to\Bbb R_+$:
$$S_k(\frak w)(j,\ell):=S_k(w_\ell)(j).$$

  Before starting the construction,
  we need some notions of {\tt block normalization.}
\subsection*{Block normalizations}
\

Suppose that $h\in\Bbb N \ \text{and}\ w\in\Bbb R_+^h$ is a block.  
\

Write  $$|h|:=h,\ M(w):=\max_{1\le j\le h}w_j,\ \ \Si(w):=\sum_{1\le j\le h}w_j\ \ \text{and}\ \ \ E(w):=\frac{\Si(w)}{|w|}$$
Note that\
$$E(w)=\int_{[1,h]\cap\Bbb N}wdP_{[1,h]\cap\Bbb N}.$$
\

The block $w\in\Bbb R_+^h$ is $\e$-{\it normalized} if
\begin{align*}S_k(w)=kE(w)(1\pm\e)\ \ \forall\ k\ge \frac{\e \Si(w)}{M(w)}. 
\end{align*}

We call the block array $\frak w\subset\Bbb R_+^h$ $\e$-{\it normalized} if each block $w \in\frak w$ is $\e$-normalized.

\subsection*{Block array distributions}
\

Let $X$ be a metric space. We'll identify the collection $\mathcal P(X)$ of Borel probabilities on $X$ with
$$\text{\tt RV}\,(X):=\{\text{\tt random variables with values in $X$}\}$$ by
$$Y\in \text{\tt RV}\,(X)\ \leftrightarrow\ \text{\tt dist.}\,(Y)\in\mathcal P(X)$$ where
$$\text{\tt dist.}\,(Y):=P\circ Y^{-1}\in\mathcal P(X)$$ in case $Y$ is defined on the probability space $(\Om,\mathcal F,P)$.

\

A {\it symmetric representation} of $Y\in\text{\tt RV}\,(X)$ is   an ordered pair $(\Om,f)$ where 
$\Om$ is a finite set and  $f:\Om\to X$ is so that
$$\text{\tt Prob}\,(Y=x)=\frac1{|\Om|}\#\,\{\om\in\Om:\ f(\om)=x\}\ \ \forall\ x\in X.$$
Evidently, the random variable $Y\in\text{\tt RV}\,(X)$ has a symmetric representation iff $Y$ is {\it rational} in the sense that there is a finite set $V\subset X$ so that
$Y\in V$ a.s. and 
$$\text{\tt Prob}\,(Y=x)\in\Bbb Q_+\ \forall\ x\in F.$$
\

Let $Y\in\text{\tt RV}\,(\Bbb R_+)$ be rational.  A $Y$-{\it distributed, $h$-block array} is a $h$-block array of form
$$\frak w\subset\Bbb R_+^h$$  with respect to which,  block averaging is a symmetric representation for $c\cdot Y$ for some $c=c(\frak w)\in\Bbb R_+$.
\

Specifically, 
$$\text{\tt Prob}\,(c\cdot Y=x)=\frac1{|\frak w|}\#\,\{w\in\frak w:\ E(w)=c\cdot x\}\ \ \forall\ x\in\Bbb R_+.$$

\

\sms{\bf Definition:\ Relative $Y$-distribution}
\

Let $Y\in\text{\tt RV}\,(\Bbb R_+)$ be rational, let  $\D>\mathcal E>0$, $h,Q\in\Bbb N$ and let
 $\frak w\subset\Bbb R_+^h\ \text{and}\ \frak w'\subset\Bbb R_+^{Qh}$ be  $Y$-distributed block arrays with 
 \
 
 $\frak w'$ refining $\frak w$, $\frak w$ $\D$-normalized and   $\frak w'$ $\mathcal E$-normalized.
 \
 
 We'll say that the pair $(\frak w,\frak w')$ is {\it relatively, $Y-(\D,\mathcal E)$-distributed} if  
 \sms (i) $m([F_{\frak w'}\ne F_{\frak w}])<\D$,
 \sms (ii) $\exists\ c(\frak w)=\g(h)\le \g(h+1)\le\dots\le\g(h')=c(\frak w')$ and $\D\ge \e_h>\e_{h+1}>\dots>\e_{Qh}=\mathcal E$
 so that $\g(k+1)-\g(k)\le\D$ and
 $$\frak u\left(\tfrac{S_k(\frak w')}{k\g(k)},Y\right)\ \ <\ \ \e_k\ \ \forall\ h\le k\le Qh.$$

 \
 
 The  proof of theorem 1 for rational random variables is based on the:
 \proclaim{Step function extension lemma}\ \ Let $Y\in\text{\tt RV}\,(\Bbb R_+)$ be rational,  let  
 $\D>0\ \text{and}\ h\in\Bbb N$.
 If $\frak w\subset\Bbb R_+^h$ is  a  $\D$-normalized,
 $Y$-distributed block array, then
\

for any  $0<\mathcal E<\D$ and  $Q\in\Bbb N$ large enough, there is a 
homogeneous $Qh$-block array $\frak w'$ refining $\frak w$ transitively so that
\

$F_{\frak w'}\ge F_{\frak w}$ and so that 
$(\frak w,\frak w')$ is  relatively $Y-(\D,\mathcal E)$-distributed.\endproclaim

 \section*{\S4  Proof of theorem 1 in the rational case}\ \ \ 
\

We first prove this case of theorem 1 assuming the  step function  extension lemma.

\

 Fix $Y\in\text{\tt RV}\,(\Bbb R_+)$. Given $\D_n\downarrow\ 0$,  with $\D_1<\frac19\min\,Y$, we build
 using the step function extension lemma iteratively, 
 a refining sequence of   block arrays $(\frak w_n)_{n\ge 1}$ with each refinement transitive and each 
$(\frak w_n,\frak w_{n+1})$ is  relatively, $Y-(\D_n,\D_{n+1})$-distributed. This gives an {\tt ESP} with distributional limit $Y$ establishing
(\dsrailways) as on page \pageref{choochoo}. 
\

To see (\Bicycle) as on page \pageref{bike}, we note that by the extension lemma, for $|\frak w|\le k\le|\frak w_{n+1}|$, we have
a coupling of 
$$\frac{S_k(\frak w_{n+1})}{k\g(k)}\ \ \text{and}\ \ Y$$
so that
$$\tfrac{S_k(\frak w_{n+1})}{k\g(k)}\ge Y-\tfrac19\min\,Y\ge\tfrac{8}9 Y$$
By monotonicity,
$$\tfrac{S_k(\frak w_{\nu})}{k\g(k)}\ge \tfrac{8}9 Y\ \ \forall\ \nu\ge n+1$$
whence
$$\tfrac{S_k(f)}{k\g(k)}\ge \tfrac{8}9 Y$$
where $F_{\frak w_\nu}\to f$ a.s..
Thus 
$$P([\tfrac{S_k(f)}{k\g(k)}<t])\le P(Y\le\tfrac{9}9 t)\ \forall\ t>0.\ \ \ \CheckedBox\ \text{\Bicycle}\ \ \qed$$

\

The rest of this  section is a proof of the  step function  extension lemma.

 \
 
 The proof is via block concatenation and perturbation.
 
 \

\proclaim{Basic lemma I}

Let $0<\D<1$ and let $w\in\Bbb R_+^h$ be $\D$-normalized. For each $$0\le \kappa\le \D E(w),\ \ \ \d>0\ \ \text{and}\ q>\frac1\D,$$ then for
$\mu\in\Bbb N$ large enough: if
 $m:=\mu q\ \text{and}\ w'\in\Bbb R_+^{mh}$ is defined  by
$$w'=w^{(\mu)}:=w^{\odot m}+\kappa qh1_{[1,mh]\cap qh\Bbb Z},$$ 
 then

 \begin{align*}&\tag{i}\text{ $w'$ is $\d$-normalized};\\ &\tag{ii}E(w')=E(w)+\kappa;\\ &\tag{iii}P(S_J(w')\ne S_J(w^{\odot m}))\le \frac{J}{qh}\ \forall\ 1\le J\le qh;
\\ &\tag{iii'}P(S_k(w')=S_k(w^{\odot m})\ \ \forall\ 1\le k\le \sqrt\D qh )\ge 1-\sqrt\D ; \\
&\tag{iv}S_k(w')=kE(w)(1\pm 2\sqrt\D) \ \forall\ \sqrt\D qh\le k\le qh;\\
&\tag{v}S_k(w')=k(E(w)+\kappa)(1\pm (\D\wedge\tfrac1k+\tfrac{\D qh}k))\ \forall\ k>qh.
\end{align*}\endproclaim
\sms{\bf Remarks.}
\

(a) Note that $F_{w'}\ge F_w$.
\

(b)  There is no contradiction between (iv) and (v) for $k\sim qh$ as the error in (iv) is at least $\frac{\kappa}{E(w)}$ which is the increment in (v).

\

\

\

\demo{Proof for $\kappa>0$}
\demo{ Proof of (i)}
\

Let $v\in\Bbb R_+^H$ be a block. We claim that 
\begin{align*}
 \tag{\Coffeecup} \frac{S_k(v)}{kE(v)}\xrightarrow[k\to\infty]{}\ 1.
\end{align*}

To see this, let $k=JH+r$ where $J\ge 1\ \text{and}\ 0\le r<H$, then
$$S_k(v)=S_{JH}(v)\pm HM(v)=JE(v)\pm HM(v)=kE(v)\pm 2HM(v)$$
whence
\begin{align*}\frac{S_k(v)}{kE(v)}&=1\pm \frac{2HM(v)}{kE(v)}\\ &\xrightarrow[k\to\infty]{}\ 1. 
\end{align*}
We have,
$$w'=w^{(\mu)}:=(w^{\pprime})^{\odot\mu}$$
where
$$w^{\pprime}:=w^{\odot q}+\kappa qh1_{\{qh\}}.$$
It follows that 
$$E(w^{(\mu)})=E(w^{\pprime})\ \text{and}\ \ M(w^{(\mu)})=M(w^{\pprime}).$$
By (\Coffeecup), $\d$-normalization  of $w'$ is obtained by  enlarging $\mu$.\ \ \ \Checkedbox

\demo{Proof of (ii)} \ \ We have

$$S_k(w')(\nu)=S_k(w^{\odot m})(\nu)+ \kappa qh\#([\nu,\nu+k-1]\cap qh\Bbb Z)\ \ \forall\ \nu\in [1,mh].$$
                                                                                                            
Therefore\
\begin{align*}
&S_{Jqh}(w')=Jq\Si(w)+J\kappa qh,\ \Si(w')=m\Si(w)+\mu\kappa qh \ \text{and}\ \ E(w')=E(w)+\kappa.\ \ \text{\Checkedbox  (ii)} \end{align*}
Also
$$ S_k(w')\le S_k(w^{\odot m})+\kappa qh\lcl\frac{k}{qh}\rcl\le S_k(w^{\odot m})+k\kappa(1+\frac{qh}k);$$ and
$$ S_k(w')\ge S_k(w^{\odot m})+\kappa qh\lfl\frac{k}{qh}\rfl\ge S_k(w^{\odot m})+k\kappa(1-\frac{qh}k).$$

\demo{ Proof of (iii) and\ (iii')}
\

$$S_k(w')=S_k(w^{\odot m})\ \ \text{\tt on}\ \ [1,mh]\setminus\bigcup_{1\le J\le\frac{m}q}(Jhq-k,Jhq]\ \ \therefore$$
$$P(S_K(w')\ne  S_K(w^{\odot m})\le\ \frac{K}{qh};\ \ \text{\Checkedbox\ (iii)}\ \ \text{and}\ $$
$$P(S_k(w')= S_k(w^{\odot m})\ \forall\ 1\le k\le \sqrt\D  qh )\ge  1- \sqrt\D .\ \ \text{\Checkedbox\ (iii')}$$
\demo{ Proof of (iv) and (v)}
\

We begin with an estimate of $S_k(w^{\odot m})$ for $k\ge \D h$.
\

\begin{align*}\tag{\S}S_k(w^{\odot m})=kE(w)(1\pm\D\wedge\tfrac{h}k)\ \ \forall\ k\ge \D h.
\end{align*}
\demo{\rm Proof of (\S)}
\

For $\D h\le k\le h,$ we have $\D\wedge\tfrac{h}k=\D$ and (\S) follows from the $\D$-normalization of $w$.
\

Let $h\le k$, then $k=Jh+r$ with $J\ge 1\ \text{and}\ r<h$ and
\begin{align*}S_k(w^{\odot m})(\nu)&=JhE(w)+\sum_{i=\nu+Jh}^{\nu+Jh+r-1}w_i\\ &=
kE(w) -rE(h)+\sum_{i=\nu+Jh}^{\nu+Jh+r-1}w_i\\ &=:kE(w)+\mathcal E.
\end{align*}
Thus
$$-\Si(w)<-rE(h)\le\mathcal E\le S_r(w)(\nu\ \mod h)\le\Si(w)$$
and
$$\frac{|\mathcal E|}{kE(w)}\le\frac{\Si(w)}{kE(w)}=\frac{h}k.$$
To see the other estimation, we use the $\D$-normalization of $w$.

\

If $r\le \frac{\D hE(w)}{ M(w)}$, then
$$|\mathcal E|\le Mr\le \D hE(w);$$
and if $r> \frac{\D hE(w)}{ M(w)}$, then by $\D$-normalization of $w$,
$$\mathcal E=-rE(w)+S_r(\nu+Jh)=-rE(w)+rE(w)(1\pm\D)=\pm\D E(w).\ \ \CheckedBox\S$$

We have 
\

\begin{align*}S_k(w')(\nu)-S_k(w^{\odot m})(\nu)=
\kappa qh\#([\nu,\nu+k-1]\cap qh\Bbb Z).\end{align*}
For $\sqrt\D qh \le k<qh,\ \#([\nu,\nu+k-1]\cap qh\Bbb Z)=0,1$
\begin{align*}S_k(w')-S_k(w^{\odot m})\le \kappa qh\le \D E(w)qh<\sqrt\D\cdot kE(w)\end{align*}
 and by (\S),
 $$S_k(w')=kE(w)(1\pm (\D\wedge\tfrac{h}k+\sqrt\D))=kE(w)(1\pm 2\sqrt\D).\ \ \CheckedBox\text{\rm (iv)} $$
\

For $k\ge qh$,
\begin{align*}S_k(w')(\nu)-S_k(w^{\odot m})(\nu)&=
\kappa qh\#([\nu,\nu+k-1]\cap qh\Bbb Z)\\ &=
\kappa qh(\frac{k}{qh}\pm 1)\\ &=\kappa k\pm \kappa qh.\end{align*}
Therefore\begin{align*}S_k(w')&=S_k(w^{\odot m})+\kappa k\pm \kappa qh\\ &=
kE(w)(1\pm \D\wedge\tfrac{h}k)+\kappa k\pm \kappa qh\\ &=
k(E(w)+\kappa)(1\pm  (\D\wedge\tfrac{h}k+ \tfrac{\kappa qh}{kE(w)}))\\ &=
k(E(w)+\kappa)(1\pm  (\D\wedge\tfrac{h}k+ \tfrac{\D qh}{k})).\ \ \CheckedBox\text{\rm (v)}
 \end{align*}
\

This proves the basic lemma.\ \ \qed

\subsection*{Example 1\ \ Constant limit random variable}\ \ \ 
\

To see how the basic lemma works, we build a  sequence of (trivial) block arrays $(\frak w_n)_{n\ge 1}$ with each
$\frak w_n=\{w^{(n)}\}$ a single block. This will give $Y\equiv 1$ as distributional limit.
\

We'll define $f^{(n)}:=w^{(n)}:\Bbb Z_{b_n}\to\Bbb R_+$ where $b_n=|w^{(n)}|$.
\

Suppose that each block $w^{(n)}$ is constructed from $w^{(n-1)}$ using the basic lemma with parameters
$$\D_n,\ \kappa_n,\ q_n,\ \mu_n,\ m_n,\ \d_n=\D_{n+1}.$$
\begin{align*}
 \tag*{\P1}\exists\ \lim_{n\to\infty}f^{(n)}=:f\in\Bbb R_+\ \ \text{a.s.}.
\end{align*}

\pf
$$P([w^{(n)}\ne w^{(n-1)}])=\frac1{q_n|w^{(n-1)}|}.$$
Since  $\sum_{n=1}^\infty\frac1{q_n|w^{(n-1)}|}<\infty$, $\exists\ N:\Om\to\Bbb N$ so that a.s., $f^{(k)}\equiv f^{(N)}\ \forall\ k\ge N$.\ \ \Checkedbox

\Par2 If $\sum_{n=1}^\infty\kappa_n=\infty$, then as $n\to\infty$,
$$E(w^{(n)})\ \sim\ \sum_{k=1}^n\kappa_k.$$

Now let $(\Om,\mathcal F,P,T)$ be the corresponding odometer and let $f:=\lim_{n\to\infty}f^{(n)}:\Om\to\Bbb R_+$.
\

Define $b:\Bbb N\to\Bbb R_+$ by
$$b(N):=NE(w^{(n)})\ \ \ \text{for}\ \ |w^{(n-1)}|<N\le |w^{(n)}|,\ n\ge 1.$$
\Par3 If $\kappa_n\to 0\ \text{and}\ \sum_{n=1}^\infty\kappa_n=\infty$, then 
$$\frac{b(n)}n\uparrow\infty,\ \ \frac{b(2n)}{b(n)}\xrightarrow[n\to\infty]{}\ 2$$
and 
\begin{align*}
\frac1{b(n)}\sum_{k=0}^{n-1}f\circ T^k\xrightarrow[n\to\infty]{\d}\ 1.
\end{align*}

\

In  example 1, the normalizing constants were directly determined by the sequence $(E(w^{(n)}))_{n\ge 1}$ of block expectations, which increased slowly.

For more complicated limit random variables (e.g. $Y\in\text{\tt RV}\,(\Bbb R_+)$ given by $P(Y=1)=P(Y=2)=\frac12$) this is no longer the case as the distributions of
the block expectations need to be considered.  A more elaborate construction procedure is necessary.

We'll need the following simultaneous version of Basic Lemma I which is an immediate consequence of it.

\

\proclaim{Basic Lemma II}
\

Let $\frak w\subset\Bbb R_+^h$ be a  $\D$-normalized $h$-block array and let $\kappa:\frak w\to \Bbb R_+$ satisfy $ 0\le \kappa(w)\le \D E(w)$.
\

For each $\d>0\ \ \text{and}\ q>\frac1\D,$ and
$\mu\in\Bbb N$ large enough: if
 $m:=\mu q$ and the   $mh$-bock array $\frak w':=\{v(w)\in\Bbb R_+^{mh}:\ w\in\frak w\}$ is defined  by
$$v(w)=w^{(\mu)}:=w^{\odot m}+\kappa(w) qh1_{[1,mh]\cap qh\Bbb Z},\ \ (w\in\frak w)$$ 
 then $\frak w'\succ\frak w\ \text{and}\ F_{\frak w'}\ge F_\frak w$ and for $w\in\frak w$,

 \begin{align*}&\tag{i}\text{ $v(w)$ is $\d$-normalized};\\ &\tag{ii}E(v(w))=E(w)+\kappa(w);\\ &\tag{iii}P(S_J(v(w))\ne S_J(w^{\odot m}))\le \frac{J}{qh}\ \forall\ 1\le J\le qh;
\\ &\tag{iii'}P(S_k(v(w))=S_k(w^{\odot m})\ \ \forall\ 1\le k\le \sqrt\D qh )\ge 1-\sqrt\D ; \\
&\tag{iv}S_k(v(w))=kE(w)(1\pm 2\sqrt\D) \ \forall\ \sqrt\D qh\le k\le qh;\\
&\tag{v}S_k(v(w))=k(E(w)+\kappa(w))(1\pm (\D\wedge\tfrac1k+\tfrac{\D qh}k))\ \forall\ k>qh.
\end{align*}\endproclaim

The next lemma is an iteration of the  procedure in Basic Lemma II to achieve larger, but  gradual  changes of the block averages $E(w)$. 
We'll use it to prove both the step function extension lemma and the step function straightening lemma.
\proclaim{Compound  lemma }
\ 

Let $0<\D<1, h\in\Bbb N$ and let $\frak w\subset\Bbb R_+^h$ be a  $\D$-normalized $h$-block array. 
Let  $\frak t:\frak w\to (1,\infty)$, then
$\forall\ \ \b>0\ \ \text{and}\ \ \ \mathcal E>0,$ and $Q\in\Bbb N$ large enough,
there is  an $\mathcal E$-normalized, $Qh$--block array
$$\frak w':=\{v(w):\ w\in\frak w\}\subset\Bbb R_+^{Qh},$$ 
numbers 
$$\d_k\ge\d_{k+1},\ \d_{Qh}<\mathcal E\ \ \text{and}\ \ 0=p_h<p_{h+1}<\dots<p_{Qh}=1,\  0\le p_{k+1}-p_k\le\b$$

so that $\frak w'\succ\frak w\ \text{and}\ F_{\frak w'}\ge F_\frak w$ 
 for each $w\in\frak w$,
\begin{align*}&\tag{ii}E(v(w))=\frak t(w) E(w);\\ 
&\tag{iii}P(S_k(v(w))=S_k(w^{\odot Q})\ \ \forall\ 1\le k\le \sqrt\D h)>1-2\sqrt\D\\ &\tag{iv}\forall\ k>\D h,\ 
 S_k(v(w))\ge kE(w)((1-p_k)+p_k\frak t(w))(1-\d_k ) \ \ \text{and}\\ &
 P([S_k(v(w))=kE(w)((1-p_k)+p_k\frak t(w))(1\pm \d_k )])\ge 1-\d_k .\end{align*}
 \endproclaim

\demo{Proof of the step function extension lemma}

Suppose that that $Y\in\text{\tt RV}\,(\Bbb R_+)$ is rational.  Let :
\bul  $(\Om,f)$ be a symmetric representation of $Y$ with $|\Om|\ge 2$, 
\bul $\frak w=\{w^{(\om)}:\ s\in\Om\}\subset\Bbb R_+^h$ be a  $\D$-normalized  block array, where $\D>0\ \text{and}\ h\in\Bbb N$ so that
$$E(w^{(\om)})=c\cdot f(\om)\ \ (\om\in\Om)$$ where $c=c(\frak w)>0$.

\

Fix  $0<\mathcal E<\D$. We'll construct for any  $Q\in\Bbb N$ large enough,  a

 $Qh$-block array $\frak w'=\{w^{\prime(s)}:\ s\in\Om\}\subset\Bbb R_+^{Qh}$ so that
 $$E(w^{\prime(s)})=c'\cdot f(s)\ \ (\om\in\Om)$$ where $c'=c(\frak w')>c(\frak w)$;
  $\frak w'\succ\frak w$ is a transitive, homogeneous extension and
$(\frak w,\frak w')$ is  relatively, $Y-(\D,\mathcal E)$-distributed.
\

The  construction is via auxiliary, intermediary block arrays $\frak w^{(1)},\frak w^{(2)},\dots,\frak w^{(N)}$ where  $N>\frac1{\mathcal E}$ is arbitrary and fixed.

\

Let $V\subset\Bbb R_+$ be the value set of $Y$ and let 
$$K>\frac{2\max\,V}{\min\,V}\ \ \text{and}\ \ N':=2(|\Om|-1)N.$$
We have that $\min_{s,t}\tfrac{Kf(t)}{f(s)}>1$ and so,
using the compound lemma, we can find $J_1>1$ and for each $s,t\in \Om$ find $\mathcal E$-normalized $w^{(s,t)}(1),\in\Bbb R_+^{J_1h}$ so that
\begin{align*}&\tag{o}E(w^{(s,t)}(1))=Kcf(t)=\tfrac{Kf(t)}{f(s)}E(w^{(s)});
\\ &\tag{i}\ P(S_k(w^{(s,t)}(1))= S_k(w^{(s)\odot J_1})\ \forall\ 1\le k\le \D J_1h)>1-\D;\\
&\tag{ii}\ c=\g(k_0)\le  \g(k_0+1)\le \dots\ \le \g(qh)=Kc;
\\ &\tag{iii} P([S_k(w^{(s,s)}(1))=k\g(k)f(s)(1\pm\D)])\ge 1-\D\ \ \forall\ k>k_0. 
\end{align*}
Here
$\g(k)=E(w^{(s)})((1-p_k)+p_kK)$ is as in the compound lemma with $\frak t\equiv K$.

\

The first intermediary block array is
$$\frak w^{(1)}=\{w^{(s,s)}(1,k):\ 1\le k\le |\Om|(N'-|\Om|+1),\ s\in \Om\}\cup\{w^{(u,v)}(1):\ u,\ v\in \Om,\ u\ne v\}$$
where $w^{(s,s)}(1,k)\ \ (1\le k\le N-1)$ is a copy of $w^{(s,s)}(1)$.
\

Next, find $J_2\ge 1$ and for each $s,t,u\in \Om,\ s\ne t$ find $w^{(s,t,u)}(2)\in\Bbb R_+^{J_2J_1h}$ so that
\begin{align*}&\text{(iii')}\ \ \ E(w^{(s,t,u)}(\nu))=cK^2f(u)=\tfrac{Kf(u)}{f(t)}E(w^{(s,t)}(1)),
\\ &\text{(iv)}\ \ \ P(S_k(w^{(s,t,u)}(2))=S_k(w^{(s,t)}(1)^{\odot J_2})\ \forall\ 1\le k\le \D J_2J_1h)>1-\D;\\
&\text{(v)}\ \ \ \ Kc=\g(k_0)\le  \g(k_0+1)\le \dots\ \le \g(qh)=K^2c
\\ &\text{(vi)}\ \ \  P([S_k(w^{(s,t,t)}(\nu))=k\g(k)f(t)(1\pm\D)])\ge 1-\D\ \ \forall\ k>k_0. 
\end{align*}
The second intermediary block array is
\begin{align*}
 &\frak w^{(2)}=\\ &\{w^{(s,s,s)}(2,k):\ 1\le k\le |\Om|(N'-2(|\Om|-1)),\ s\in \Om\}\cup\{w^{(s,t,t)}(2),w^{(s,s,t)}(2):\ s,t\in \Om,\ s\ne t\}
\end{align*}

where $w^{(s,s,s)}(2,k)\ \ (1\le k\le N-2)$ is a copy of $w^{(s,s,s)}(2)$.

\

Recurse this, to find $J_2,J_3,\dots,J_N$ and for each $2\le\nu\le N,\ s_1,s_2,\dots,s_\nu\in \Om$,  $w^{(s_1,s_2,\dots,s_\nu)}(\nu)\in\Bbb R_+^{h^{(\nu-1)}}$  where
$h^{(\nu)}:=hJ_1J_2\dots J_\nu$; so that
\begin{align*}&\text{(iii')}\ \ \ E(w^{(s_1,s_2,\dots,s_\nu)}(\nu))=cK^\nu f(s_\nu)=\tfrac{Kf(s_\nu)}{f(s_{\nu-1})}E(w^{(s_1,s_2,\dots,s_{\nu-1})}(\nu-1)),\\ &
\text{(iv)}\ \ \ P(S_k(w^{(s_1,s_2,\dots,s_\nu)}(\nu))=S_k(w^{(s_1,s_2,\dots,s_{\nu-1})}(\nu-1))^{\odot J_\nu})\ \forall\ 1\le k\le \D h^{(\nu)})>1-\D;\\
&\text{(v)}\ \ \ \ K^{\nu-1}c=\g(k_0)\le  \g(k_0+1)\le \dots\ \le \g(qh)=K^{\nu}c;
\\ &\text{(vi)}\ \ \  P([S_k(w^{(s_1,s_2,\dots,s_{\nu-2},t,t)}(\nu))=f(t)k\g(k)(1\pm\D)])\ge 1-\D\ \ \forall\ k>k_0. 
\end{align*}

The $\nu^{\text{\tt\tiny th}}$ intermediary block array is
$$\frak w^{(\nu)}=\{w^{(s^\nu)}(\nu,k):\ 1\le k\le |\Om|(N'-\nu(|\Om|-1)),\ s\in \Om\}\cup\bigcup_{j=1}^{\nu-1}\{w^{(s^j,t^{\nu-j})}(\nu):\ s,t\in \Om,\ s\ne t\}$$
where $w^{(s^\nu)}(\nu,k)\ \ (1\le k\le N-\nu)$ is a copy of $w^{(s^\nu)}(\nu)$.

In particular,
$$\frak w^{(N)}=\{w^{(s^N)}(N,k):\ 1\le k\le |\Om|(N'-N(|\Om|-1)),\ s\in \Om\}\cup\bigcup_{j=1}^{N-1}\{w^{(s^j,t^{N-j})}(N):\ s,t\in \Om,\ s\ne t\}$$
where $w^{(s^N)}(N,k)\ \ (1\le k\le N-N)$ is a copy of $w^{(s^N)}(N)$.

\

 Now set $\frak w'=\{w^{'(s)}:\ s\in \Om\}$ where
 
 $$w^{'(s)}:=\(\bigodot_{k=1}^{N(|\Om|-1)} w^{(s^N)}(N,k)\ \odot\  \bigodot_{t\in \Om\setminus\{s\}}\bigodot_{j=1}^N w^{(t^{N-j},s^j)}(N)\)^{\odot T}$$
 where $T$ is chosen large enough to ensure $\mathcal E$-normalization.
 \

This is as advertised.\ \ \Checkedbox

\section*{\S5 General case of Theorem 1 and Theorem 2}

We now complete the proof of theorem 1 by constructing an {\tt ESP} with an arbitrary  $Y\in\text{\tt RV}\,(\Bbb R_+)$ as distributional limit. 

\

For this, we need to approximate an arbitrary $Y\in\text{\tt RV}\,(\Bbb R_+)$ with rational random variables in a controlled manner.
 \subsection*{Splittings}
 \

 A {\it splitting} of the finite set $\Om$ is a surjection $\pi:\Xi\to\Om$  defined on another  finite set $\Xi$ 
so that  $P_{\Om}=P_\Xi\circ\pi^{-1}$. 
\

Equivalently, $\#\pi^{-1}\{x\}=\frac{\#\Xi}{\#\Om}\ \forall\ x\in\Om$.
 
 \

\

 Let  the compact metric space $([0,\infty],\rho)$ be as before, let $\pi:\Xi\to\Om$ be a splitting  and let $(\Om,f)$,  $(\Xi,g)$
 be symmetric representations.
 \
 
 We'll say, for $\e>0$, that   
 $(\Xi,g)$   $\e$-{\it splits}  $(\Om,f)$ via  $\pi:\Xi\to\Om$  
 if 
 $$E_\Xi(\rho(g,f\circ\pi)):=\frac1{\#\Xi}\sum_{u\in\Xi}\rho(g(u),f(\pi(u))))\ \ <\ \ \e$$
 and we'll call  $\pi:\Xi\to\Om$  the (associated) {\it $\e$-splitting}.
 
 \

 Note that if $Z$ has a symmetric representation which 
 $\e$-splits some symmetric representation of $Y$, then
 $\frak v(Y,Z)<\e.$

 \proclaim{Splitting approximation lemma}
 \
 
 Let $Y\in\text{\tt RV}\,(\Bbb R_+)$,
 \ \ then  
 $\forall\ \ \e_k\downarrow 0$ there is  a 
 sequence $(Y_1,Y_2,\dots)$ of rational random variables  on $\Bbb R_+$ with a  nested sequence of symmetric representations $(\Om_k,f_k)$ 
 so that 
 \sms {\rm (o)}   $\frak v(Y_k,Y)<\e_k\ \forall\ k\ge 1$;
 \sms {\rm (i)}  $(\Om_{k+1},f_{k+1})$ $\e_k$-splits $(\Om_{k},f_{k})\ \forall\ k\ge 1$.
 \sms {\rm (ii)} $\exists\ R>0$ so that  $P_{\Om_k}(Y_k<t)\le  \text{\tt Prob}\,(Y<t)\ \forall\ t\in (0,R),\ k\ge 1$.
 \endproclaim\demo{Proof}\ \ \ Considering $Y$ as a random variable on the compact metric space $([0,\infty],\rho)$, we  
 let $\mu:=\text{\tt dist}\,(Y)\in \mathcal P([0,\infty])$. There is a non-decreasing  map  $\Phi:[0,1]\to [0,\infty]$ so that
 $\mu=\l\circ\Phi^{-1}$ where $\l$ is Lebesgue measure on $[0,1]$. Let $\G\subset [0,1]$ be the collection of discontinuity points of $\Phi$.
 By monotonicity, this set is at most countable.
 \

 Let $Z:=\{0,1\}^\Bbb N$ equipped with the product, discrete topology, and let $B:Z\to [0,1]$ be the ``binary expansion map''
 $$B((x_1,x_2,\dots)):=\sum_{k=1}^\infty\frac{x_k}{2^k}.$$
 It follows that the collection of discontinuity points of   $\Psi:=\Phi\circ B:Z\to [0,\infty]$ is $\widetilde{\G}=B^{-1}\G$.
 This set is also at most countable.
 \
 
 We have 
 $$\mu=\nu\circ\Psi^{-1}$$ where
 $\nu=\prod(\tfrac12,\tfrac12)\in\mathcal P(Z)$. 
 
 \
 
 By the above,
 $$\Phi(\sum_{k=1}^{n-1}\frac{x_k}{2^k}+\frac1{2^n})\xrightarrow[n\to\infty]{}\Psi(x_1,x_2,\dots)\ \text{for $\nu$-a.e.}\ \ (x_1,x_2,\dots)\in Z$$
 (indeed $\forall\ (x_1,x_2,\dots)\notin\widetilde{\G}$).

 Now, for $n\ge 1$, let $Z_n:=\{0,1\}^n$, define $\psi_n:Z_n\to [0,1]$ by $$\psi_n(x_1,x_2,\dots,x_n):=\Phi(\sum_{k=1}^{n-1}\frac{x_k}{2^k}+\frac1{2^n}).$$ 
We have that for $\nu$-a.e. $(x_1,x_2,\dots)\in Z$,
 $$\psi_n(x_1,x_2,\dots,x_n)\xrightarrow[n\to\infty]{}\Psi(x_1,x_2,\dots).$$
 Define the restriction maps $\pi_n:Z\to Z_n\ \text{and}\ \pi_n^{n+m}:Z_{n+m}\to Z_n$ by 
 $$\pi_n(x_1,x_2,\dots)=(x_1,x_2,\dots,x_n)\ \text{and}\ \pi_n^{n+m}(x_1,x_2,\dots,x{n+m})=(x_1,x_2,\dots,x_n),$$
 then $\pi_n^{n+m}:Z_{n+m}\to Z_n$ is a splitting and 
by Egorov's theorem, along a sufficiently sparse subsequence $n_k\uparrow\infty$, we have
$$\int_Z\rho(\psi_{n_k}\circ\pi_{n_k},\Psi)d\nu<\frac{\e_k}2$$ whence
$$E_{Z_{n_{k+1}}}(\rho(\psi_{n_k}\circ\pi_{n_k}^{n_{k+1}},\psi_{n_{k+1}}))<\e_k.$$
                                                                                 
Thus 
 $$\Om_k:=Z_{n_k},\ f_k:=\psi_{n_k}\ \text{and}\ \text{\tt dist}\,(Y_k):=P_{\Om_k}\circ f_k^{-1}\in\mathcal P(\Bbb R_+)$$
 are as required for (i), which entails (o).
 \
 
To see (ii) we note that 
$$\psi_n(x_1,x_2,\dots,x_n)\ge \Psi((x_1,x_2,\dots)$$
whenever $(x_1,x_2,\dots,x_n)\ne\mathbb{1}$. 
Let 
$$R:=\Phi(\sum_{j=1}^{n_1-1}\frac1{2^j})=\Phi(1-\frac1{2^{n_1}})\le \Phi(\sum_{j=1}^{n_k-1}\frac1{2^j})\ \forall\ k\ge 1.$$
If $k\ge 1\ \text{and}\ \psi_{n_k}(x_1,x_2,\dots,x_{n_k})<R$ then $(x_1,x_2,\dots,x_{n_k})\ne\mathbb{1}$ and
$\psi_{n_k}(x_1,x_2,\dots,x_{n_k})\ge \Psi((x_1,x_2,\dots)$.

\

Since $f_k=\psi_{n_k}$, for $t\in (0,R)$
$$P_{\Om_k}([f_k\le t])\le \nu([\Psi\le t])=P(Y\le t).\ \ \ \CheckedBox\text{(ii)}$$

  \proclaim{Step function straightening lemma}
 \
 
 Let $Y,\ Z\in\text{\tt RV}\,(\Bbb R_+)$ be rational with symmetric representations $(\Om,f)$ and  $(\Xi,g)$ respectively.
 \

Suppose that $\mathcal E,\ \D>0$ and that $(\Xi,g)$ $\mathcal E$-splits $(\Om,f)$ with  $\mathcal E$-splitting $\Phi:\Xi\to\Om$.
\

Let $\frak w=\{w(\om):\ \om\in\Om\}\subset\Bbb R_+^h$
be a $\D$-normalized,  $Y$-distributed, $h$-block array with $E(w(\om))=c(\frak w)f(\om)\ \forall\ \om\in\Om$.
 \

Then for each $Q\in\Bbb N$ large enough and $\eta>0$, $\exists$
 a \ $\mathcal E$-normalized,  $(\Xi,g)$-distributed, $Qh$-block array
 $$\frak b=\{b(\xi):\ \xi\in\Xi\}\subset\Bbb R_+^{Qh},$$ so that  
 $$F_{\frak b}\ge F_{\frak w}\ \ \ \text{and}\ \ \ m([F_{\frak b}\ne F_{\frak w}])<\mathcal E,$$
 and
$$\b(h)\le \b(h+1)\le\dots\le\b(Qh),\ \ \ \b(k+1)-\b(k)\le\eta,$$ $$0=q_h<q_{h+1}<\dots<q_{Qh}=1,\ \d_h\ge\d_{k+1}\ge\dots\ge\d_{Qh},\  \d_{Qh}<\mathcal E$$
so that for $h\le k\le Qh$,
\begin{align*}&
 S_k(b(\xi))\ge k\b(k)((1-q_k)f(\Phi(\xi))+q_kg(\xi))(1- \d_k )\\ &
  P([S_k(b(\xi))=k\b(k)((1-q_k)f(\Phi(\xi))+q_kg(\xi))(1\pm \d_k )])\ge 1-\d_k\\ &
  \frak v(\tfrac{S_k(\frak b)}{k\b(k)},Z)<\mathcal E+\D.\end{align*}
\endproclaim
\demo{Proof}
 \
 
 Let $\Phi:\Xi\to\Om$ be so that
 $$P_\Xi\circ\Phi^{-1}=P_\Om\ \text{and}\ \ E_\Xi(\rho(f\circ\Phi,g))<\mathcal E.$$
 For $\xi\in\Xi$, let $v(\xi):=w(\Phi(\xi))\in\frak w$ and consider the block array 
 $$\widetilde{\frak w}:=\{v(\xi):\ \xi\in\Xi\}.$$ 
 
 Note that $E(v(\xi))=cf(\Phi(\xi))$. In order to use the  compound lemma, define $\frak t:\Xi\to (1,\infty)$ by
                                                                                                           
 $$\frak t(\xi):=\frac{Kg(\xi)}{f(\Phi(\xi))}\ \ \text{where}\ \ K>\max_{\xi\in\Xi}\frac{f(\Phi(\xi))}{g(\xi)}$$
 so that $\frak t>1$.
 \
 
 By the compound lemma
 for $Q\ge 1$ large enough, 
there is  an $\mathcal E$-normalized, $Qh$--block array
$$\frak b=\{b(\xi):\ \xi\in\Xi\}\subset\Bbb R_+^{Qh},$$ 
numbers 
$$\d_k\ge\d_{k+1},\ \d_{Qh}<\mathcal E\ \ \text{and}\ \ 0=p_h<p_{h+1}<\dots<p_{Qh}=1,\ p_{k+1}-p_k<\eta$$

so that
 for each $\xi\in\Xi$,
\begin{align*}
&E(b(\xi))=\frak t(\xi) E(v(\xi))=c(\frak w)f(\Phi(\xi));\\ 
&P(S_k(b(\xi))=S_k(v(\xi)^{\odot Q})\ \ \forall\ 1\le k\le \D h)>1-2\D\end{align*}
and $\forall\ k>\D h$,
\begin{align*}
 P([S_k(b(\xi))=kE(b(\xi))((1-p_k)+p_k\frak t(\xi))(1\pm \d_k )])\ge 1-\d_k.\end{align*}
Next, for $\xi\in\Xi$,
 $$E(b(\xi))((1-p_k)+p_k\frak t(\xi))=c(\frak w)(1-p_k)f(\Phi(\xi))+Kp_kg(\xi)).$$
 Let 
 $$\b(k):=c(\frak W)(p_k+(1-p_k)K),\ q_k:=\frac{Kp_k}{p_k+(1-p_k)K},$$
 then
 $$0=q_h<q_{h+1}<\dots<q_{Qh}=1$$ and
 $$E(b(\xi))((1-p_k)+p_k\frak t(\xi))=\b(k)((1-q_k)f(\Phi(\xi))+q_kg(\xi)).$$
 
 Thus, with probability $\ge 1-\d_k$,
 $$\rho(\frac{S_k(b(\xi))}{k\g(k)},(1-q_k)f(\Phi(\xi))+q_kg(\xi))<\d_k$$
and
 \begin{align*}E_\Xi(\rho(\frac{S_k(b(\xi))}{k\g(k)},g(\xi))&\le 2\d_k+E_\Xi(\rho(f\circ\Phi,g))\\ &\le\d_k+\mathcal E.
 \end{align*}
The inequality $F_{\frak b}\ge F_{\frak w}$ follows from monotonicity.\ \ \Checkedbox  

\

 \demo{Proof of theorem 1}
 \
 
 Fix $\e_n\downarrow 0,\ \sum_{n=1}^\infty\e_n<\infty$ and
use the splitting approximation lemma to obtain
 a  sequence $(Y_1,Y_2,\dots)$ of rational random variables on $\Bbb R_+$  
  with a  nested sequence of symmetric representations $(\Om_k,f_k)$ 
 so that 
 \sms {\rm (o)}   $\frak v(Y_k,Y)<\e_k\ \forall\ k\ge 1$;
 \sms {\rm (i)}  $(\Om_{k+1},f_{k+1})$ $\e_k$-splits $(\Om_{k},f_{k})\ \forall\ k\ge 1$.
 \sms {\rm (ii)} $\exists\ R>0$ so that  $P_{\Om_k}(Y_k<t)\le  \text{\tt Prob}\,(Y<t)\ \forall\ t\in (0,R),\ k\ge 1$.
\

Using the step function extension- and straightening lemmas (respectively), we next, construct   sequences 
$(\frak v_n)_n\ \text{and}\ (\frak e_n)_n$ of $Y_n$-distributed $h_n$- and $k_n$-block arrays (respectively) so that
$$\frak v_n\ \prec\ \frak w_n\ \ \prec \frak v_{n+1}\ \ \text{and}\ \ F_{\frak v_n}\le F_{\frak w_n}\le F_{\frak v_{n+1}}$$ 
and a slowly varying sequence $(\g(k))_k,\ \g(k+1)-\g(k)\to 0$
so that with $b(k):=k\g(k)$, for some $r>0$
\sms (iii) $m([F_{\frak v_n}\ne F_{\frak w_n}])<\e_n\ \text{and}\ m([F_{\frak w_n}\ne F_{\frak v_{n+1}}]) <\e_{n+1};$
\sms (iv) $\tfrac{S_k(\frak w_n)(\xi)}{b(k)}\ge rf_n(\xi)\forall\ h_n<k\le h_{n+1}$ where $\frak w_n=\{w(\xi):\ \xi\in\Om_n\}$,
\sms (v)  $\frak v(\tfrac{S_k(\frak w_n)}{b(k)},g)<\e_n\ \ \forall\ h_n<k\le h_{n+1}.$

\

Let $$\xbmt:=\varprojlim_{n\to\infty}\frak W_n\ \text{and}\ f:=\lim_{n\to\infty}F_{\frak W_n,\frak w_n},$$ then $(X,\B,m,T,f)$ is an  
{\tt ESP} with distributional limit $Y$.

\

Moreover, if $h_n<k\le h_{n+1}$, and $t\in (0,R)$ then $S_k(f)\ge S_k(F_{\frak w_n})$ whence
$$[S_k(f)\le tb(k)]\subset [S_k(F_{\frak w_n})\le tb(k)]$$ whence by (iv),
$$P([S_k(f)\le tb(k)])\le P([\tfrac{S_k(\frak w_n)(\xi)}{b(k)}\le t])\le P(Y_n\le \frac{t}r)\le P(Y\le \frac{t}r).\ \ \CheckedBox$$

\

\demo{ Proof of theorem 2}
\

We use the odometer construction of theorem 1 to prove theorem 2.
\

Let $Y\in\text{\tt RV}\,(\Bbb R_+)$ and let $(\Om,\mathcal F,P,\tau)$ be an {\tt EPPT}. We must exhibit a measurable function $\phi: \Om\to\Bbb R_+$ so that the
{\tt ESP} $(\Om,\mathcal F,P,\tau,\phi)$ has distributional limit $Y$. 
\

Now fix as above, an odometer $\xbmt$ with $f:X\to\Bbb R_+$ measurable so that $(X,\B,m,T,f)$ satisfies (\dsrailways)  in theorem 1 (on page \pageref{choochoo})
with distributional limit $Y$ and
$1$-regularly varying normalizing constants $b(n)_{n\ge 1}$.
\

 \ By the odometer factor proposition, there is a set $\Om_0\in\mathcal F,\ P(\Om_0)>0$ so that the induced {\tt EPPT} 
$(\Om_0,\mathcal F\cap\Om_0,P_{\Om_0},\tau_{\Om_0})$ 
has $\xbmt$ as a factor.

Let  $\phi:(\Om_0,\mathcal F\cap\Om_0,P_{\Om_0},\tau_{\Om_0})\to\xbmt$ be the factor map and define $\pi:\Om\to\Bbb R$ by
$$\phi=f\circ\pi\ \text{on}\ \Om_0\ \ \text{and}\ \ \phi\equiv 0\ \ \text{off}\ \Om_0.$$
We have that
\begin{align*}
\frac1{b(n)}\sum_{k=0}^{n-1}\phi\circ \tau_{\Om_0}^k\xrightarrow[n\to\infty]{P_{\Om_0}-\frak d}\ Y.
\end{align*}
Now let $\kappa:\Om_0\to\Bbb N$ be the first return time of $\tau$ to $\Om_0$ and let $\kappa_n:=\sum_{j=0}^{n-1}\kappa\circ\tau_{\Om_0}^j$
(the $n^{\text{\tiny th}}$ return time of $\tau$ to $\Om_0$), then on $\Om_0$, 
$$\sum_{k=0}^{n-1}\phi\circ \tau_{\Om_0}^k\equiv \sum_{j=0}^{\kappa_n-1}\phi\circ \tau^j.$$
By Birkhoff's theorem, $\kappa_n\sim\frac{n}{P(\Om_0)}$ a.s. on $\Om_0$ and so by monotonicity and $1$-regular variation of $b(n))_{n\ge 1}$,
\begin{align*}
\frac1{b(n)}\sum_{k=0}^{n-1}\phi\circ \tau^k\xrightarrow[n\to\infty]{P_{\Om_0}-\frak d}\ P(\Om_0)Y
\end{align*}
whence by Eagleson's theorem,
\begin{align*}
\frac1{b(n)}\sum_{k=0}^{n-1}\phi\circ \tau^k\xrightarrow[n\to\infty]{\frak d}\ P(\Om_0)Y.\ \ \ \CheckedBox
\end{align*}
\section*{\S6  New examples in infinite ergodic theory}

We begin by reviewing:
\subsection*{Kakutani skyscrapers and inversion}
\

As in \cite{Kakutani}, the {\it  skyscraper} over the $\Bbb N$-valued {\tt SP} $(\Om,\mathcal F,P,S,f)$ is the {\tt MPT} $\xbmt$ defined by
$$X=\{(x,n):\ x\in \Om,\ 1\le n\le f(x)\},$$
$$\B=\sigma\{A\x \{n\}:n\in\Bbb N,\
A\in\mathcal F\cap [f\ge n]\},\ m(A\x\{n\})=P(A),$$ and
$$T(x,n)=\begin{cases} &(Sx,f)\text{ if }n=f(x),\\
&(x,n+1)\text{ if }1\le n\le f(x)-1.\end{cases}$$
The skyscraper {\tt MPT} is always conservative as $\bigcup_{n\ge 1}T^{-n}\Om\x\{1\}=X$ and its ergodicity is equivalent to that of
$(\Om,\mathcal F,P,S)$. Any invertible {\tt CEMPT} $\xbmt$ is isomorphic to the skyscraper over a {\it first return time} {\tt SP} 
$(\Om,\mathcal B\cap\Om ,m_\Om,T_\Om,\varphi_\Om)$ where $\varphi_\Om(x):=\min\,\{n\ge 1:\ T^nx\in\Om\}$ is the {\it  first return time} 
which is finite for a.e. $x\in\Om$ by conservativity,
$T_\Om(x):=T^{\varphi_\Om(x)}$ is the {\it induced transformation} on $\Om$ which is a {\tt PPT}.\

\

Let $\xbmt$ be an invertible {\tt CEMPT} let $\Om\in\B,\ m(\Om)=1$ and consider the {\tt return time  stochastic process} on $\Om$:
{\sms  $(\Om,\mathcal B\cap\Om ,m_\Om,T_\Om,\varphi_\Om)$ where $\varphi_\Om(x):=\min\,\{n\ge 1:\ T^nx\in\Om\}$.}
\

Distributional limits with regularly varying normalizing constants are transferred between the return time {\tt SP} and the Kakutani skyscraper 
by means of the following
 \proclaim{ Inversion proposition\ \ \cite{distlim}}
 \
 
 \ \ \ Let  $a(n)$ be $\g$-regularly varying with $\g\in (0,1]\ \text{and}$ fix $\Om\in\mathcal F$, then 
for $Y$  a rv on $(0,\infty)$:
$$\tfrac1{a(n)}S_n(1_\Om)\overset{\mathfrak d}\lra Ym(\Om)\ \ \iff\ 
\tfrac{\varphi_{n}}{a^{-1}(n)}\overset{\mathfrak d}\lra (\tfrac1{m(\Om)Y})^{\frac1\g}$$
 where $\varphi_n=\sum_{k=0}^{n-1}\varphi_\Om\circ T_\Om^k$.\endproclaim

 \demo{Proof of Theorem 3}\ \ Fix $Y\in\text{\tt RV}\,(\Bbb R_+)$,  let $(\Om,\mathcal F,P,S,f)$ be a $\Bbb N$-valued {\tt ESP} and let
 $b(n)$ be $1$-regularly varying so that
 \begin{align*}
 &\frac1{b(n)}\sum_{k=0}^{n-1}f\circ T^k\xrightarrow[n\to\infty]{\frak d}\ \frac1Y\\ &
 P([\sum_{k=0}^{n-1}f\circ T^k<xb(n)])\le P(\tfrac1Y\le t)\ \ \forall\ t>0\ \text{small}\ \text{and}\ \ n\ge 1\ \text{large}.
\end{align*}
 These exist by theorem 1. Now let $\xbmt$ be the Kakutani skyscraper over
 $(\Om,\mathcal F,P,S,f)$. By inversion,
 \begin{align*}
 &\tag{\Football}\frac{S_n^{(T)}}{b^{-1}(n)}\xyr[n\to\infty]{\frak d}\ Y\ \ \text{and} \\ &\tag{\dsaeronautical}
 m_\Om([S_n^{(T)}(1_\Om)>xb^{-1}(n)])\le P(Y\text{ $\ge x$})\ \ \forall\ y>1,\ n\ge 1\ \ \text{large}.\ \ \ \CheckedBox
\end{align*}
\
 
 \subsection*{Rational ergodicity properties}
 \
 
Now let $\a>0$ and let $K\subset\Bbb N$ be a subsequence. 
\

We'll say that the {\tt CEMPT} $\xbmt$ is {\it $\a$-rationally ergodic along $K$} if for some $\Om\in\B,\ 0<m(\Om)<\infty$, we have
\begin{align*}\tag{{\tt $\a$-RE}$_K$}\int_A\left(\frac{S_n(1_B)}{a(n)}\right)^\a dm\xrightarrow[n\to\infty,\ n\in K]{}\ m(A)m(B)^\a \ \forall\ A,\ B\in\B(\Om)\end{align*}
where $a(n)=a_{\a,\Om}(n):=\frac1{m(\Om)^{1+\frac1\a}}(\int_\Om(S_n(1_\Om)^\a dm)^\frac1\a$.
\

We'll say that $\xbmt$ is {\it $\a$-rationally ergodic} if it is  $\a$-rationally ergodic along $\Bbb N$ and 
{\it subsequence  $\a$-rationally ergodic} if it is  $\a$-rationally ergodic along some $K\subset\Bbb N$.

Properties like this have been considered in   \cite{AdamsSilva} and \cite{Roblin}.
\

Standard techniques show that 
$\Om\in\B,\ 0<m(\Om)<\infty$ satisfies ({\tt $\a$-RE}$_K$)  iff
$$\left\{\left(\frac{S_n(1_\Om)}{a_{\a,\Om}(n)}\right)^\a:\ n\in K\right\}$$ is uniformly integrable on $\Om$, and, if nonempty, the collection
$$R_{\a,K}(T):=\{\Om\in\B:\  0<m(B)<\infty\  \text{\tt  satisfying  ({\tt $\a$-RE}$_K$)}\}$$
is a dense $T$-invariant hereditary ring. 
\

Moreover
$a_{\a,\Om}(n)\sim a_{\a,\Om'}(n)$ along $K$ whenever  $\Om,\ \Om'\in R_{\a,K}(T)$
\

{ We'll call the {\tt CEMPT} $\xbmt$  {\it  $\infty$-rationally ergodic along $K$} if for some $\Om\in\B,\ 0<m(\Om)<\infty$, we have
\begin{align*}\tag{{\tt BRE}$_K$}\sup_{n\in K}\left\|\frac{S_n(1_\Om)}{a_{1,\Om}(n)}\right\|_{L^\infty(\Om)}<\infty.\end{align*}
Analogously to as above, if nonempty, the collection
$$R_{\infty,K}(T):=\{\Om\in\B:\  0<m(B)<\infty\  \text{\tt  satisfying  ({\tt BRE}$_K$)}\}$$
is a dense $T$-invariant hereditary ring. It is contained in $R_{\a,K}(T)\ \forall\ \a>0$. 
\

The condition  {\tt $\infty$-rational ergodicity  along $\Bbb N$} is aka  {\tt bounded rational ergodicity}. For more information and examples, see   \cite{BRE}.}
\subsection*{$\a$-return sequence}\ \ 
We define the $\a$-{\it return sequence} of an $\a$-rationally ergodic {\tt CEMPT} $\xbmt$ as the growth rate
$$a_{n,\a}(T)\sim a_{\a,\Om}(n)\ \ \Om\in R_\a(T).$$

\

It is also possible to define  ``subsequence $\a$-return sequence'' for a subsequence $\a$-rationally ergodic {\tt CEMPT}.

\

Note that 
\bul $1$-rational ergodicity is equivalent to weak rational ergodicity as in \cite{RE} with  $R_1(T)=R(T)$ and $a_{n,1}(T)\sim a_n(T)$;

\bul $2$-rational ergodicity implies rational ergodicity;
\bul for { $0<\a\le\infty$}, $\a$-rational ergodicity
 implies $\b$-rational ergodicity for each $\b\in (0,\a)$;
\bul pointwise dual ergodic transformations are  $\a$-rationally ergodic $\forall\ 0<\a<\infty$ (this follows from the existence of moment sets).

\

Let $\xbmt$ be distributionally stable with limit $Y\in\text{\tt RV}\,(\Bbb R_+)$.
\

 { For $\a\in\Bbb R_+$, set $\|Y\|_\a:=E(Y^\a)^{\frac1\a}\le\infty$ and $$\|Y\|_\infty:=\sup\,\{t>0:\ P(Y>t)>0=\lim_{\a\to\infty}\|Y\|_\a\le\infty.$$
\bul For $0<\a\le\infty$,  if $T$ is $\a$-rationally ergodic, then $\|Y\|_\a<\infty$ and if $\a\in\Bbb R_+$, then 
$a_{n,\a}(T)\sim \|Y\|_\a a_{n,Y}(T)$. 
\bul If $\|Y\|_\a=\infty$, then  $T$ is not subsequence,  $\a$-rationally ergodic.}

\subsection*{Example:\ distributional stability $\nRightarrow$\ $\a$-rational ergodicity}
\

Let $Y\in\text{\tt RV}\,(\Bbb R_+)$ be so that  $E(Y^\a)=\infty\ \forall\ \a>0$. By  theorem 3, there is a distributionally stable {\tt CEMPT} 
$(X,\B,m,T)$ with ergodic limit $Y$
with $a_{n,Y}(T)$ $1$-regularly varying.  By the above $\forall\ \a>0,\ \ T$ is not subsequence,  $\a$-rationally ergodic.

\

For a given {\tt CEMPT} $\xbmt$, we consider the collection
$$I(T):=\{\a>0:\ T\ \text{is $\a$-rationally ergodic}\}.$$
It follows from the above that $I(T)$ must be an interval, either empty, or $\Bbb R$, or of  form $(0,a)$ or $(0,a]$ for some $a\in (0,\infty]$.

\

We conclude this paper by showing that all possibilities occur.

\

\proclaim{Lemma}
\

\ \  Let $\xbmt$ be distributionally stable with ergodic limit $Y\in\text{\tt RV}\,(\Bbb R_+)$ and $a_{n,Y}(T)$ 
$1$-regularly varying. Suppose that $\Om\in\B,\ m(\Om)=1$ satisfies (\dsaeronautical) as on page \pageref{aero}, then
$T$ is $\a$-rationally ergodic iff $\|Y\|_\a<\infty$ and in this case, when $\a<\infty$, $a_{n,\a}(T)\sim E(Y^\a)^\frac1\a a_{n,Y}(T)$.\endproclaim
\

\demo{Proof of $\|Y\|_\a<\infty\ \Lra$  $\a$-{\tt RE}}
\

\ { We only consider the case $0<\a<\infty$. The case where $\a=\infty$ is  easy.}  We claim first that
$$\{\Phi_n:=(\tfrac{S_n(1_\Om)}{a_{n,Y}(T)})^\a:\ n\ge 1\}$$ is a uniformly integrable family in $L^1(\Om)$.
\

Now, since $E(Y^\a)<\infty$,  we have by monotone convergence and Fubini's theorem that
$$\rho(t):=\int_t^\infty P(Y^\a>s)ds=E(1_{[Y^\a>t]}Y^\a)\xyr[t\to\infty]{}\ 0.$$
By (\dsaeronautical) (page \pageref{aero}), 
\begin{align*}\int_\Om 1_{[\Phi_n>t]}\Phi_ndm &=\int_t^\infty m([\Phi_n>s])ds\\ &\le
28\int_t^\infty P(Y^\a>s)ds\\ &=:\rho(t)
\end{align*}
whence
$$\sup_{n\ge 1}\int_\Om 1_{[\Phi_n>t]}\Phi_ndm\le \rho(t)\xyr[t\to\infty]{}\ 0$$
and the family is uniformly integrable.
\

Next by (\Football) as on page \pageref{football}, for $A,B\in\B(\Om)\ \text{and}\ x>0$,

$$\int_A(\tfrac{S_n(1_B)}{a_{n,Y}(T)})^\a\wedge xdm \xyr[n\to\infty]{}\ m(A)E((m(B)Y)^\a\wedge x).$$
\

Moreover, $E(m(B)Y)^\a\wedge x) \xyr[x\to\infty]{}\ m(B)^\a E(Y^\a)$. To estimate the error,
\begin{align*}0&\le \int_A(\tfrac{S_n(1_B)}{a_{n,Y}(T)})^\a dm -\int_A(\tfrac{S_n(1_B)}{a_{n,Y}(T)})^\a\wedge xdm\\ &\le 
 \int_A(\tfrac{S_n(1_B)}{a_{n,Y}(T)})^\a 1_{[(\frac{S_n(1_B)}{a_{n,Y}(T)})^\a>x]}dm\\ &\le
 \int_\Om 1_{[\Phi_n>x]}\Phi_ndm\\ &\le \ \rho(x)\xrightarrow[x\to\infty]{}\ 0.
\end{align*}
Standard arguments now show that
$$\int_A(\tfrac{S_n(1_B)}{a_{n,Y}(T)})^\a dm\xyr[n\to\infty]{}\ m(A)m(B)^\a E(Y^\a).\ \ \CheckedBox$$
{ Note that a boundedly rationally ergodic transformation $T$ has $I(T)=(0,\infty]$ and a pointwise, 
dual ergodic transformation $T$ with return sequence which is regularly varying with index $\g<1$
has as ergodic limit a $\g$-Mittag-Leffler random variable (see \cite{distlim}) which is unbounded but has moments of all orders, whence $I(T)=(0,\infty)$.}
\

The following completes the picture (and is also a strengthening of \cite{AdamsSilva}):
\proclaim{Proposition}\ \   For each $a\in \Bbb R_+$ there are  distributionally stable {\tt MPT}s $T_o\ \text{and}\ T_c$ with
 $I(T_o)=(0,a)$   or $I(T_c)=(0,a]$. 
\endproclaim 

 \demo{Proof of the Proposition}\ \  To construct $T_o$ with $I(T_o)=(0,\a)$ fix a 
 $Y\in\text{\tt RV}\,(\Bbb R_+)$ so that $E(Y^t)<\infty\ \forall\ t<\a$ but $E(Y^\a)=\infty$ and construct $T$ as in the theorem 3.
 
 \

To construct $T_c$ with $I(T_c)=(0,\a]$ the same but using a 
 $Z\in\text{\tt RV}\,(\Bbb R_+)$ so that $E(Z^\a)<\infty$ but $E(Z^t)=\infty\ \forall\ t>\a$.\ \Checkedbox
 \


\begin{bibdiv}
\begin{biblist}

\bib{RE}{article}{
      author={Aaronson, Jon},
       title={Rational ergodicity and a metric invariant for {M}arkov shifts},
        date={1977},
        ISSN={0021-2172},
     journal={Israel J. Math.},
      volume={27},
      number={2},
       pages={93\ndash 123},
      review={\MR{0584018}},
}

{\bib{BRE}{article}{
      author={Aaronson, Jon.},
       title={Rational ergodicity, bounded rational ergodicity and some
  continuous measures on the circle},
        date={1979},
        ISSN={0021-2172},
     journal={Israel J. Math.},
      volume={33},
      number={3-4},
       pages={181\ndash 197 (1980)},
         url={http://dx.doi.org/10.1007/BF02762160},
        note={A collection of invited papers on ergodic theory},
      review={\MR{571529 (81f:28012)}},
}}

\bib{distlim}{article}{
      author={Aaronson, Jon},
       title={The asymptotic distributional behaviour of transformations
  preserving infinite measures},
        date={1981},
        ISSN={0021-7670},
     journal={J. Analyse Math.},
      volume={39},
       pages={203\ndash 234},
         url={http://dx.doi.org/10.1007/BF02803336},
      review={\MR{632462 (82m:28030)}},
}

\bib{IET}{book}{
      author={Aaronson, Jon},
       title={An introduction to infinite ergodic theory},
      series={Mathematical Surveys and Monographs},
   publisher={American Mathematical Society},
     address={Providence, RI},
        date={1997},
      volume={50},
        ISBN={0-8218-0494-4},
      review={\MR{1450400 (99d:28025)}},
}

\bib{AS}{article}{
      author={Aaronson, Jon},
      author={Sarig, Omri},
       title={Exponential chi-squared distributions in infinite ergodic
  theory},
        date={2014},
        ISSN={0143-3857},
     journal={Ergodic Theory Dynam. Systems},
      volume={34},
      number={3},
       pages={705\ndash 724},
         url={http://dx.doi.org/10.1017/etds.2012.160},
      review={\MR{3199789}},
}

\bib{gendistlim}{article}{
      author={Aaronson, Jon},
      author={Weiss, Benjamin},
       title={Generic distributional limits for measure preserving
  transformations},
        date={1984},
        ISSN={0021-2172},
     journal={Israel J. Math.},
      volume={47},
      number={2-3},
       pages={251\ndash 259},
         url={http://dx.doi.org/10.1007/BF02760521},
      review={\MR{738173 (85e:28024)}},
}

\bib{AZ}{article}{
      author={Aaronson, Jon},
      author={Zweim{\"u}ller, Roland},
       title={Limit theory for some positive stationary processes with infinite
  mean},
        date={2014},
        ISSN={0246-0203},
     journal={Ann. Inst. Henri Poincar\'e Probab. Stat.},
      volume={50},
      number={1},
       pages={256\ndash 284},
         url={http://dx.doi.org/10.1214/12-AIHP513},
      review={\MR{3161531}},
}

\bib{AdamsSilva}{article}{
      author={{Adams}, T.~M.},
      author={{Silva}, C.~E.},
       title={{Weak Rational Ergodicity Does Not Imply Rational Ergodicity}},
        date={2015-02},
     journal={ArXiv e-prints},
      eprint={1502.06566},
}

\bib{Aldous}{article}{
      author={Aldous, D.~J.},
      author={Eagleson, G.~K.},
       title={On mixing and stability of limit theorems},
        date={1978},
     journal={Ann. Probability},
      volume={6},
      number={2},
       pages={325\ndash 331},
      review={\MR{0517416}},
}

\bib{ADDS}{article}{
      author={Avila, A.},
      author={Dolgopyat, D.},
      author={Duryev, E.},
      author={Sarig, O.},
       title={The visits to zero of a random walk driven by an irrational
  rotation},
        date={2015},
        ISSN={0021-2172},
     journal={Israel J. Math.},
      volume={207},
      number={2},
       pages={653\ndash 717},
         url={http://dx.doi.org/10.1007/s11856-015-1186-4},
      review={\MR{3359714}},
}

\bib{Billingsley}{book}{
      author={Billingsley, Patrick},
       title={Convergence of probability measures},
     edition={Second},
      series={Wiley Series in Probability and Statistics: Probability and
  Statistics},
   publisher={John Wiley \& Sons, Inc., New York},
        date={1999},
        ISBN={0-471-19745-9},
         url={http://dx.doi.org/10.1002/9780470316962},
        note={A Wiley-Interscience Publication},
      review={\MR{1700749}},
}

\bib{BGT}{book}{
      author={Bingham, N.~H.},
      author={Goldie, C.~M.},
      author={Teugels, J.~L.},
       title={Regular variation},
      series={Encyclopedia of Mathematics and its Applications},
   publisher={Cambridge University Press},
     address={Cambridge},
        date={1987},
      volume={27},
        ISBN={0-521-30787-2},
      review={\MR{898871 (88i:26004)}},
}

\bib{Bromberg}{article}{
      author={{Bromberg}, M.},
       title={{Ergodic Properties of the Random Walk Adic Transformation over
  the Beta Transformation}},
        date={2015-11},
     journal={ArXiv e-prints},
      eprint={1511.02482},
}

\bib{RB-MD}{article}{
      author={Burton, Robert},
      author={Denker, Manfred},
       title={On the central limit theorem for dynamical systems},
        date={1987},
        ISSN={0002-9947},
     journal={Trans. Amer. Math. Soc.},
      volume={302},
      number={2},
       pages={715\ndash 726},
         url={http://dx.doi.org/10.2307/2000864},
      review={\MR{891642 (88i:60039)}},
}

\bib{Chacon}{inproceedings}{
      author={Chacon, R.~V.},
   booktitle={Proceedings of the fifth berkeley symposium on mathematical
  statistics and probability, volume 2: Contributions to probability theory,
  part 2},
   publisher={University of California Press},
     address={Berkeley, Calif.},
       pages={335\ndash 360},
}

\bib{DK}{article}{
      author={Darling, D.~A.},
      author={Kac, M.},
       title={On occupation times for {M}arkoff processes},
        date={1957},
        ISSN={0002-9947},
     journal={Trans. Amer. Math. Soc.},
      volume={84},
       pages={444\ndash 458},
      review={\MR{0084222 (18,832a)}},
}

\bib{Eagleson}{article}{
      author={Eagleson, G.~K.},
       title={Some simple conditions for limit theorems to be mixing},
        date={1976},
        ISSN={0040-361x},
     journal={Teor. Verojatnost. i Primenen.},
      volume={21},
      number={3},
       pages={653\ndash 660},
      review={\MR{0428388 (55 \#1409)}},
}

\bib{Feller}{book}{
      author={Feller, William},
       title={An introduction to probability theory and its applications.
  {V}ol. {II}},
   publisher={John Wiley \& Sons, Inc., New York-London-Sydney},
        date={1966},
      review={\MR{0210154}},
}

\bib{Friedman}{book}{
      author={Friedman, Nathaniel~A.},
       title={Introduction to ergodic theory},
   publisher={Van Nostrand Reinhold Co., New York-Toronto, Ont.-London},
        date={1970},
        note={Van Nostrand Reinhold Mathematical Studies, No. 29},
      review={\MR{0435350 (55 \#8310)}},
}

\bib{Kakutani}{article}{
      author={Kakutani, Shizuo},
       title={Induced measure preserving transformations},
        date={1943},
     journal={Proc. Imp. Acad. Tokyo},
      volume={19},
       pages={635\ndash 641},
      review={\MR{0014222}},
}

\bib{Levy}{book}{
      author={L{\'e}vy, P.},
       title={Calcul des probabilit{\'e}s},
      series={PCMI collection},
   publisher={Gauthier-Villars},
        date={1925},
         url={https://books.google.co.il/books?id=8\_FLAAAAMAAJ},
}

\bib{Renyi}{article}{
      author={R{\'e}nyi, A.},
       title={On mixing sequences of sets},
        date={1958},
        ISSN={0001-5954},
     journal={Acta Math. Acad. Sci. Hungar.},
      volume={9},
       pages={215\ndash 228},
      review={\MR{0098161 (20 \#4623)}},
}

\bib{Roblin}{article}{
      author={Roblin, Thomas},
       title={Sur l'ergodicit\'e rationnelle et les propri\'et\'es ergodiques
  du flot g\'eod\'esique dans les vari\'et\'es hyperboliques},
        date={2000},
        ISSN={0143-3857},
     journal={Ergodic Theory Dynam. Systems},
      volume={20},
      number={6},
       pages={1785\ndash 1819},
         url={http://dx.doi.org/10.1017/S0143385700000997},
      review={\MR{1804958 (2001m:37061)}},
}

\bib{Shields-cutstack}{article}{
      author={Shields, Paul~C.},
       title={Cutting and stacking: a method for constructing stationary
  processes},
        date={1991},
        ISSN={0018-9448},
     journal={IEEE Trans. Inform. Theory},
      volume={37},
      number={6},
       pages={1605\ndash 1617},
         url={http://dx.doi.org/10.1109/18.104321},
      review={\MR{1134300 (92k:28031)}},
}

\bib{Shields-book}{book}{
      author={Shields, Paul~C.},
       title={The ergodic theory of discrete sample paths},
      series={Graduate Studies in Mathematics},
   publisher={American Mathematical Society, Providence, RI},
        date={1996},
      volume={13},
        ISBN={0-8218-0477-4},
         url={http://dx.doi.org/10.1090/gsm/013},
      review={\MR{1400225 (98g:28029)}},
}

\bib{Skorohod}{article}{
      author={Skorohod, A.~V.},
       title={Limit theorems for stochastic processes},
        date={1956},
        ISSN={0040-361x},
     journal={Teor. Veroyatnost. i Primenen.},
      volume={1},
       pages={289\ndash 319},
      review={\MR{0084897}},
}

\bib{TZ}{article}{
      author={Thaler, M.},
      author={Zweim{\"u}ller, R.},
       title={Distributional limit theorems in infinite ergodic theory},
        date={2006},
        ISSN={0178-8051},
     journal={Probab. Theory Related Fields},
      volume={135},
      number={1},
       pages={15\ndash 52},
         url={http://dx.doi.org/10.1007/s00440-005-0454-3},
      review={\MR{2214150 (2007e:60017)}},
}

\bib{JPT-BW}{article}{
      author={Thouvenot, Jean-Paul},
      author={Weiss, Benjamin},
       title={Limit laws for ergodic processes},

        date={2012},
        ISSN={0219-4937},
     journal={Stoch. Dyn.},
      volume={12},
      number={1},
       pages={1150012, 9},
         url={http://dx.doi.org/10.1142/S0219493712003596},
      review={\MR{2887924}},
}

\bib{Vas}{article}{
      author={Vasershtein, L.~N.},
       title={Markov processes over denumerable products of spaces describing
  large system of automata},
        date={1969},
     journal={Problemy Pereda\v ci Informacii},
      volume={5},
      number={3},
       pages={64\ndash 72},
      review={\MR{0314115 (47 \#2667)}},
}

\bib{Volny}{article}{
      author={Voln{\'y}, Dalibor},
       title={Invariance principles and {G}aussian approximation for strictly
  stationary processes},
        date={1999},
        ISSN={0002-9947},
     journal={Trans. Amer. Math. Soc.},
      volume={351},
      number={8},
       pages={3351\ndash 3371},
         url={http://dx.doi.org/10.1090/S0002-9947-99-02401-0},
      review={\MR{1624218 (99m:60062)}},
}

\end{biblist}
\end{bibdiv}
\end{document}